\documentclass[reqno]{amsart}

\usepackage{capt-of}
\usepackage[skip=6pt,font=footnotesize]{caption}
\usepackage{amsmath}
\usepackage{amsfonts}
\usepackage{amssymb}
\usepackage{graphicx}
\usepackage{algorithmic}
\usepackage{todonotes}
\usepackage{xcolor}
\usepackage{array}
\usepackage{extarrows}
\usepackage{wasysym}
\usepackage{tikz}
\usepackage{pgfplots}
\usetikzlibrary{plotmarks}
\usepackage{pgfplotstable}
\pgfplotsset{compat=1.12}
\usetikzlibrary{external}

\graphicspath{ {./images/} }
\usepackage{mathrsfs}
\usepackage{subfig}

\numberwithin{equation}{section}

\usepackage{hyperref}
\hypersetup{colorlinks=true,allcolors=black}

\usepackage{amsmath}

\usepackage{longtable}
\usepackage{stmaryrd}
\usepackage{mathrsfs}
\usepackage{xcolor}

\usepackage{bm}         

\usepackage{esint}	 
\usepackage{nicefrac}	
\newcommand{\minidist}{\hspace{1pt}}
\usepackage{letltxmacro}

\usepackage{mhchem} 
\newcommand{\ie}{\emph{i.e.} }

\newcommand{\tennabla}{\boldsymbol{\nabla}}

\renewcommand{\div}{\text{div}}
\newcommand{\tf}{\text{tf}}

\newcommand{\vol}[1]{\text{vol}\{#1\}}				  				   
\newcommand{\area}[1]{\text{area}\{#1\}}

\newcommand{\ustilde}[1]{\tilde{\underset{s}{#1}}\,\!}

\usepackage{mhchem} 

\newcommand{\An}{\text{An}}
\newcommand{\El}{\text{El}}
\newcommand{\Sep}{\text{Sep}}
\newcommand{\Cat}{\text{Cat}}

\newcommand{\A}{\mathtt{A}}
\newcommand{\E}{\mathtt{E}}
\newcommand{\C}{\mathtt{C}}

\newcommand{\D}{\mathtt{D}}
\renewcommand{\P}{\mathtt{P}}
\renewcommand{\S}{\mathtt{S}}
\newcommand{\vecx}{\mathbf{x}}
\newcommand{\vecy}{\mathbf{y}}
\newcommand{\vece}{\mathbf{e}}

\newcommand{\vecchi}{\boldsymbol{\chi}}
\newcommand{\tenpi}{\boldsymbol{\pi}}

\newcommand{\tenId}{\minidist \textbf{Id}}
\newcommand{\us}[1]{\underset{s}{#1}\,\!}

\newcommand{\SCL}{\text{SCL}}

\newcommand{\EA}{{\E_A}}
\newcommand{\iAE}{{\A\E}}

\newcommand{\iCD}{{\C\D}}
\newcommand{\iCE}{{\C\E}}

\newcommand{\vecn}{\mathbf{n}}
\newcommand{\EC}{{\E_C}}
\newcommand{\ES}{{\E_S}}

\newcommand{\Etot}{{\E,\text{\tiny{tot}}}}
\newcommand{\Alat}{{\A,\text{\tiny{lat}}}}
\newcommand{\ESref}{{\ES,\text{\tiny{ref}}}}
\newcommand{\ECref}{{\EC,\text{\tiny{ref}}}}
\newcommand{\Eq}{{\E_q}}

\newcommand{\Sq}{{\S_q}}
\newcommand{\Aq}{{\A_q}}
\renewcommand{\O}{\mathcal{O}}
\newcommand{\vecJ}{\mathbf{J}}       	
\newcommand{\flux}{\vecJ}

\newcommand{\Aell}{{\A_\ell}}

\newcommand{\AM}{{\A_M}}
\newcommand{\Ae}{{\A_e}}

\newcommand{\ustlambda}{\us{\tilde{\lambda}}}

\usepackage{eurosym}

\newcommand{\mathR}{\mathbb{R}}

\newcommand{\eps}{\varepsilon}

\let\phi\varphi

\newcommand{\cearrow}{\rightleftharpoons}

\newcommand{\Ell}{\mathscr{L}}
\newcommand{\Uone}{\mathfrak{1}}
\newcommand{\Um}{\text{m}}

\newcommand{\Ucm}{\text{cm}}
\newcommand{\Unm}{\text{nm}}
\newcommand{\Umum}{\mu\text{m}}

\newcommand{\UA}{\text{A}}
\newcommand{\UAh}{\text{Ah}}

\newcommand{\Ul}{\text{L}}

\newcommand{\UV}{\text{V}}

\newcommand{\Uh}{\text{h}}

\newcommand{\Umol}{\text{mol}}

\newcommand{\Us}{\,\text{s}}

\newcommand{\const}{\text{const.}}

\newcommand{\rref}{\text{\tiny{ref}}}
\newcommand{\tot}{\text{\tiny{tot}}}
\newcommand{\lat}{\text{\tiny{lat}}}

\newcommand{\units}[1]{\,/\,#1\,}
\newcommand{\uunits}[1]{\left[ #1 \right]}

\newcommand{\kT}{k_\text{\tiny{B}}T\,}

\renewcommand{\exp}[1]{\text{e}^{#1}\,}

\renewcommand{\ln}[1]{\text{ln}\left( #1 \right)}

\newcommand{\ceE}[2][\,\!]{
\ensuremath{\ce{#2}\big|_{\E^{#1}}}
}

\newcommand{\ceA}[2][\,\!]{
\ensuremath{\ce{#2}\big|_{\A}^{#1}}
}

\ifpdf
\hypersetup{
    pdftitle={A Modeling Framework for Efficient Reduced Order Simulations},
    pdfauthor={M. Landstorfer, M. Ohlberger, S. Rave and  M. Tacke}
}
\fi

\title[A Modeling Framework for Efficient Reduced Order Simulations]{A Modeling Framework for Efficient Reduced Order Simulations of Parametrized Lithium-Ion Battery Cells}

\author{M. Landstorfer$^1$} 
\address{$^1 \,$Weierstrass-Institute, Mohrenstrasse 39, 10117 Berlin, 
Germany, \textup{ \texttt{manuel.landstorfer @wias-berlin.de}}}

\author{M. Ohlberger$^2$}    
\address{$^2 \,$Center for Nonlinear Science \& Applied Mathematics Muenster, Einsteinstrasse 62, 48149 Muenster, Germany, \textup{ \texttt{\{mario.ohlberger, stephan.rave, marie.tacke\}@uni-muenster.de}}}

\author{S. Rave$^2$}
\author{M. Tacke$^2$}

\date{\today}

\thanks{This work was supported by the “Bundesministerium für Bildung und Forschung” (BMBF) through the research grants No. 05M18BCA and 05M18PMA. In addition, we acknowledge funding by the Deutsche Forschungsgemeinschaft under Germany’s Excellence Strategy EXC 2044-390685587, Mathematics Münster: Dynamics – Geometry – Structure (M. Ohlberger, S. Rave and M.Tacke) and EXC 2046: MATH+: Berlin Mathematics Research Center (AA4-8, M. Landstorfer).}


\begin{document}

\begin{abstract}
In this contribution we present a new modeling and simulation framework for parametrized Lithium-ion battery cells. We first derive a new continuum model  
for a rather general intercalation battery cell on the basis of non-equilibrium thermodynamics. In order to efficiently evaluate the resulting parameterized non-linear system of partial differential equations the reduced basis method is employed. The reduced basis method is a model order reduction technique on the basis of an incremental hierarchical approximate proper orthogonal decomposition approach and empirical operator interpolation. 
The modeling framework is particularly well suited to investigate and quantify degradation effects of battery cells. Several numerical experiments are given to demonstrate the scope and efficiency of the modeling framework.
\end{abstract}

\maketitle


\section{Introduction}
Lithium-ion batteries (LIBs) are a key component of our modern society, with applications ranging from medical devices via consumer electronics to electric vehicles and aerospace industry. The further development of LIBs is based on various aspects, namely safety, cost, storage capacity and degradation stability. This research and development is assisted by mathematical models, which are capable of simulating the complex behavior of LIBs on various degrees of spatial and temporal resolution \cite{Newman,Landstorfer_2019,Latz:2011aa}. Mathematical models based on continuum thermodynamics allow, for example, the  simulation of charging and discharging processes (cycling), yielding the cell voltage $E$  as function of the capacity $Q$ (or status of charge $y$) and the cycling rate $C_h$. This quantities are typically determined in galvanostatic electrochemical measurements, enabling a comparison between theoretical predictions and experimental data \cite{Landstorfer_2019}. 

Such models can then be used to investigate and quantify degradation effects of LIBs \cite{D1CP00359C}, which are experimentally studied with periodic charging and discharging over $N$ cycles. Experimental measurements seek to determine the dependency of the cell voltage $E$ as function of cycling rate $C_h$, cycle number $n$, and status of charge, \ie $E = E(y;C_h, n)$, to systematically quantify different ageing aspects. This requires in general a huge batch of cells and massive amounts of measuring times, e.g. cycling a cell at $C_h=1$ for $1000$ cycles requires at least $80$ days of continuous electrochemical measurements. Variations of materials or material combinations subsequently increase measuring times exponentially.  

Model based simulations can help to reduce this lab time by reducing the number of batch cells, cycle numbers and material combinations. Aging effects can be represented in a continuum model with different approaches \cite{PELLETIER2017158,HAN2019100005,BARRE2013680}, either by full spatial and temporal resolution of a specific effect, e.g. cracking due to intercalation stress, or by cycle number dependent parameters. The latter approach requires an evolution equation for parameter variations with respect to the cycle number $n$, which can itself either be upscaled from some microscopic model or determined from experimental snapshots. Once parametrized, such a model can in principle predict the cycling behavior for various cycle rates and cycle numbers, \ie $E = E(y;C_h, n)$, on the basis of numerical simulations. 

Quite similar as the experimental time expenditure increases for multiple simultaneous parameter variations, so does 
 computation time for numerical simulations. Modern mathematical tools, however, allow to reduce this computational time significantly by setting up a reduced basis (RB) method. Hence the combination of a (continuum) mathematical modeling, parametrized degradation functions, numerical simulations and reduced basis methods yield a versatile toolbox for the investigation of battery aging on the basis of electrochemical data. 
 
In this work we derive a mathematical model framework for a rather general intercalation battery cell on the basis of non-equilibrium thermodynamics \cite{dGM84}. It considers three porous phases, namely a porous intercalation anode, a porous separator phase, and a porous intercalation cathode. Each porous phase consists of an electrolyte phase, which is based on a rigorously validated material model \cite{Landstorfer2016187}, accounting for solvation effects \cite{Dreyer:2014fk}, incompressibility constraints, diffusion and conductivity. The active intercalation phase of the electrodes accounts for intercalation lithium in terms of a specific chemical potential function, and solid state diffusion with a concentration dependent diffusion coefficient \cite{Landstorfer_2019}. Furthermore, the conductive additive phase is considered, where electron transport is modeled as a simple Ohmic law. All three phases are coupled through a reaction boundary condition, where a special emphasis is put on thermodynamic consistency \cite{Landstorfer01012017,C5CP03836G}. Subsequent spatial homogenization techniques \cite{doi:10.1137/0523084} for the porous structure yield an effective, non-linear coupled partial differential equation (PDE) system for the lithium concentration in each phase, the Lithium-ion concentration in the electrolyte, and the electrostatic potential in each phase.  Building on this, we present a RB method for recurring numerical simulations of the parameterized PDE model, equipped with various degradation models expressed in terms of cycle number dependent parameters.

We discretize the resulting nonlinear system of PDEs by the finite element method in space and the backward Euler method in time. 
The computational studies to identify critical parameters for estimating the aging process of batteries require many evaluations of the finite element system with different parameter settings and thus involve a large amount of time and experimental effort.  
Therefore, we employ the RB method in order to further reduce the parameterized discretized battery model to obtain a reduced order model (ROM) that is cheap to evaluate.  For an introduction and overview on recent developments in model order reduction we refer to the monographs and collections
\cite{RB_Rozza, RB_Manzoni, MR3672144, MR3701994}. 
For time-dependent problems, the POD-Greedy method \cite{Pod_Greedy} 
defines the {\em Gold-Standard} for systems, where rigorous and cheap to evaluate a posteriori error estimates are available. As this is not the case for the non-linear battery model at hand, we employ the hierarchical approximate proper orthogonal decomposition (HAPOD) \cite{HAPOD} in this contribution.

As RB methods rely on so called efficient offline/online splitting, they need to be combined with supplementary interpolation methods in case of non-affine parameter dependence or non-linear differential equations. The empirical interpolation method (EIM) \cite{eim} and its various generalizations are key technologies with this respect. In this contribution we employ the empirical operator interpolation as introduced in \cite{eim_operator, eim_gen_operator}.

Several related model order reduction approaches have already been applied for 
battery simulation, as e.g. in \cite{White, Lass_Volkwein, Zhang_2012}, where the authors apply model reduction techniques for Newman-type Lithium-ion battery models \cite{Newman}.
Further results on model order reduction in the context of battery models, including resolved electrode geometry and Butler-Volmer kinetics can be found in \cite{Wesche_Volkwein,Stephan_Zhang, Stephan_Felix, Stephan, multibat} and \cite{Xia}.\\

The article is outlined as follows. 
In Section \ref{sect:derivation_model} we derive the entire mathematical model for a porous battery cell. 
The approximation and model order reduction of the resulting electrochemical battery model is presented in Section  \ref{sec:headings}.  
In Section \ref{sec:results} we give further details on the implementation of the modeling framework and demonstrate the scope and efficiency of the approach in several numerical experiments (Subsection \ref{sec:exp1}--\ref{sec:exp3}).
Finally, in Section \ref{sec:conclusion} we summarize the results of this contribution.

\section{Derivation of the porous electrode model}\label{sect:derivation_model}
In this section we derive the a new continuum mathematical model for a porous battery cell. 
After setting up domains and proper definitions in \ref{par:dpomains}, we state the chemical potential functions for all phases in \ref{par:variables} and briefly discuss their derivation as well as some consequence of the electro-neutrality condition. In \ref{sub:transport} we state then the corresponding transport equations for each phase in the porous electrode, where \ref{par:variables} puts a special emphasis on the intercalation reaction boundary condition and its formulation on basis of non-equilibrium surface thermodynamics. Section \ref{par:balance_equations} covers then the full PDE system of an un-homogenized porous electrode, setting the basis for spatial homogenization in section \ref{ssub:homogen}. Introducing proper scalings in \ref{eq:subsection_non-dim_1} yields a general homogenized equation framework in \ref{ssub:general_homogenization_framework} for a porous multi-phase electrode. Section \ref{ssub:homogenized_battery_model} then covers the full homogenized battery model, where proper scalings are introduced on the basis of the $1-$C current density. An overview of all parameters is given in \ref{ssub:parameters}.

\subsection{Domains, definitions and initial scaling}
\label{par:dpomains}
\begin{figure}[h] 
    \centering
    \includegraphics[width=\textwidth]{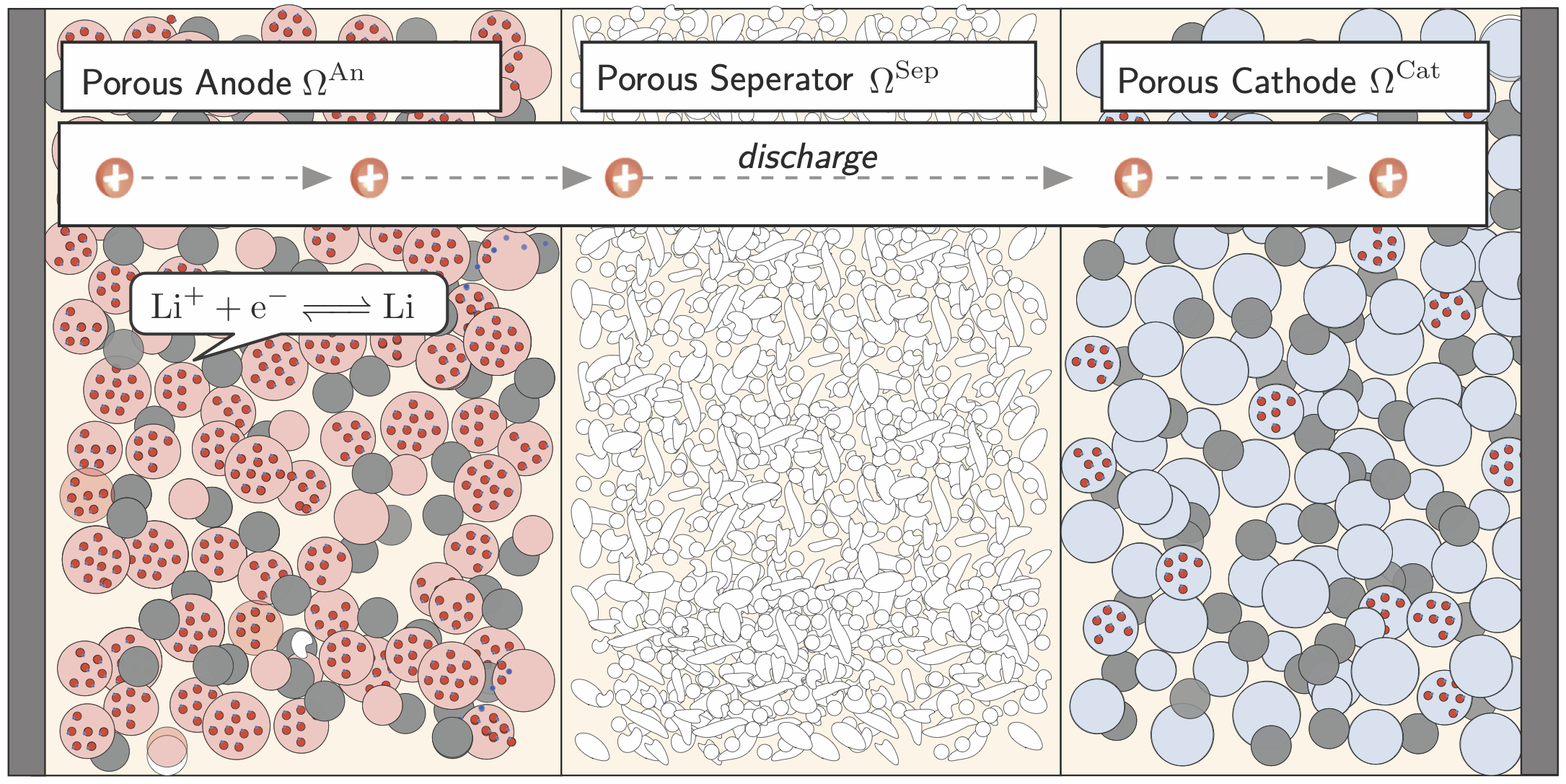}
    \caption{Sketch of the porous electrochemical unit cell. During discharge, lithium ions flow from the anode to the cathode, while electrons drive an external electrical consumer.}
        \label{fig:porous_cell_level}s
\end{figure}

We seek to model a porous electrochemical unit cell, consisting of a porous anode $\Omega^\An \subset \mathR^3$, a porous separator $\Omega^\Sep  \subset \mathR^3$, and a porous cathode $\Omega^\Cat \subset \mathR^3$ (see Fig \ref{fig:porous_cell_level}.). The intercalation electrodes $\Omega^j, j = \An,\Cat$ consist themself of three phases, namely the active particle phase $\Omega_\A^j$, the conductive additive phase $\Omega_\C^j$, and the electrolyte phase $\Omega_\E^j$, with $\Omega^j = \bigcup_{i \in \{\A,\C,\E \}} \Omega_i^j$. The union of the active phase and the conductive additive is frequently called solid phase $\Omega^j_\S = \bigcup \Omega_\A^j \cup \Omega_\C^j$. The porous separator consists of an electrolyte phase $\Omega_\E^\Sep$ and an polymeric additive $\Omega_\P^\Sep$, with $\Omega^\Sep = \Omega_\E^\Sep \cup \Omega_\P^\Sep$. The whole electrolyte phase is further denoted by $\Omega_\E = \bigcup_{j \in \{\An,\Sep,\Cat\}} \Omega_\E^j$, the active phase as $\Omega_\A = \bigcup_{j \in \{\An,\Cat\}} \Omega_\E^j$, and the solid phase $\Omega_\S = \bigcup_{j \in \{\An,\Cat\}} \Omega_\S^j$.  The interface $\Sigma_{\A,\E} = \Omega_\A \cap \Omega_\E$ captures the actual surface $\Sigma_\A$ of the active particle as well as the electrochemical double layer forming at the interface, \ie $\Sigma_{\A,\E} = \Omega_\A^\SCL \cup \Sigma_\A \cup \Omega_\E^\SCL$. The domains $\Omega_\E$ and $\Omega_\A$ are thus electro-neutral, and we refer to \cite{Newman:1965aa,Landstorfer01012017,C5CP03836G} for details on the derivation. A similar argument holds for the interface $\Sigma_\iCE = \Omega_\C \cap \Omega_\E$.

In each phase of each domain, we have balance equations, which are coupled through respective boundary conditions modeling the intercalation reaction. We seek to employ periodic homogenization to derive a homogenized PDE model for the electrochemical unit cell $\Omega = \Omega^\Cat \cup \Omega^\Sep \cup \Omega^\An$. In order to do so, we require a specific scaling with respect to the two essential length scales arising in the problem, (i) the length scale $\Ell$ of the macroscopic width $W$ of the unit cell, \ie $\Omega = \Sigma \times [0,W]$ with $\Sigma \subset \mathR^2$ being the area of the deflectors, and (ii) the length scale $\ell$ of the intercalation particle radii  $r_\A$.  This yields a small parameter $\eps = \frac{\Ell}{\ell}$ with respect to which we can perform a multi-scale asymptotic expansion and thus a periodic homogenization.  

\subsubsection{Variables, chemical potentials and parameter (equilibrium) } 
\label{par:variables}
The active particle $\Omega_\A$ is a mixture of electrons \ce{e-}, intercalated cations \ce{C} and lattice ions \ce{M^+}, and the electrolyte a mixture of solvated cations \ce{C^+}, solvated anions \ce{A-} and solvent molecules \ce{S}. The respective species densities are denoted with $n_\alpha(\vecx,t), \vecx \in \Omega_i$. Further, we denote with
\begin{align}
        \mu_\alpha = \frac{\partial \psi}{\partial n_\alpha} \qquad , \qquad  i=\A,\E,~ \alpha = \EA,\EC,\ES,{\A_C},\Ae,\AM~,
\end{align}
the chemical potential of the constituents, which are derived from a free energy density \cite{Muller:1985fk,dGM84} $\psi = \psi_\A + \psi_\E$ with $\psi_\A=\hat \psi_\A(n_\Ae,n_{\A_C},n_\AM)$ of the active particle and $\psi_\E=\hat \psi(n_\ES,n_\EA,n_\EC)$ of the electrolyte phase \cite{Muller:1985fk,Dreyer:2018aa,Landstorfer01012017,Landstorfer2016187}. \\

\paragraph{\textbf{Electrolyte}} 
\label{par:electrolyte}
For the electrolyte we consider exclusively the material model \cite{Dreyer:2014fk,Landstorfer2016187,C3CP44390F} of an incompressible liquid electrolyte accounting for solvation effects, \ie
\begin{align}
        \label{eq:electrolyte_mu}
        \mu_\alpha = g_{\alpha,\rref} + \kT \ln{y_\alpha} + v_\alpha (p-p_\rref) \quad \alpha = \ES,\EA,\EC,
\end{align}
with mole fraction
\begin{align}
        y_\alpha = \frac{n_\alpha}{n_\Etot}~,
\end{align}
molar concentration $n_\alpha$, and total molar concentration of the mixture (with respect to the number of mixing particles \cite{Landstorfer2016187})
\begin{align}
        n_\Etot = n_\ES + n_\EA + n_\EC.
\end{align}
Note that $n_\ES$ denotes the number of \emph{free} solvent molecules, whereas $n_\EA$ and $n_\EC$ are the densities of the solvated ions. This is crucial for various aspects of the thermodynamic model, and we refer to \cite{Dreyer:2014fk,Landstorfer2016187,LANDSTORFER201856,dgm2018} for details. Overall, the material model for the electrolyte corresponds to an incompressible mixture with solvation effects and we assume further 
\begin{align}
        \frac{v_\EC}{v_\ESref} = \frac{m_\EC}{m_\ES} \quad \text{and} \quad \frac{v_\EA}{v_\ES} = \frac{m_\EA}{m_\ES}
\end{align}
whereby the incompressibility constraint \cite{Dreyer:2014fk,Landstorfer2016187,C3CP44390F} implies also a conservation of mass, \ie
\begin{align}
        \label{eq:incomp_constraint}
        \sum_{\alpha} v_\alpha n_\alpha = 1 \quad \Leftrightarrow \quad \sum_{\alpha} m_\alpha n_\alpha = \rho = \frac{m_\ES}{v_\ESref} = \const. 
\end{align}
The molar volume of the solvent is related to the mole density $n_\ESref$ of the pure solvent as 
\begin{align}
        v_\ES = (n_\ESref)^{-1}~.
\end{align}
Note further that the partial molar volumes $v_\alpha$ and the molar masses $m_\alpha$ of the cation and anion are related to the solvation number $\kappa_\EC$ and $\kappa_{\A_C}$, respectively. We assume that partial molar volume of the ionic species is mainly determined by the solvation shell, which seems reasonable for large solvents like DMC in comparison to the small ions like \ce{Li+}. We proceed thus with the assumption
\begin{align}
        v_\EC = \kappa_\E \cdot v_\ES \quad \text{and} \quad v_\EA = \kappa_\E \cdot v_\ES. 
\end{align}
The electrostatic potential in the electrolyte phase is denoted as  $\phi_\E(\vecx,t)$.\\

\paragraph{\textbf{Active Phase}} 
\label{par:active_phase}
We discuss exemplarily one electrode active phase, which is in section \ref{par:balance_equations} applied to the anode and cathode phase. For the active particle phase $\Omega_\A$, we consider exemplarily a classical lattice mixture model  \cite{Cahn:1958aa,Cahn:1959ab,Garcia:2004aa,DREYER:2011aa,Dreyer:2011ab,M.-Landstorfer:2011aa,Landstorfer:2013ly,doi:10.1021/ar300145c}, which we term in the following \emph{ideal electrode} material. Note that the whole modeling as well as the numerical procedure can directly be adapted to other chemical potential functions $\mu_{\A_C} = \mu_{\A_C}(y_\A)$. The chemical potential of intercalated lithium is stated as 
\begin{align}
        \label{eq:AC_mu}
        \mu_{\A_C} = g_{\A_C} + \kT f_\A(y_\A) ~, \quad \text{with} \quad     f_\A(y_\A) &:= \ln{ \frac{y_\A}{1\! +\!  y_\A} } + \gamma_\A   \cdot (2 y_\A \!-\! 1)~,
\end{align}
with mole fraction 
        \begin{align}
                y_\A = \frac{n_{\A_C}}{n_\Alat} 
                \end{align}
of intercalated cations in the active phase. The number density $n_\Aell$ of lattice sites is constant, which corresponds to an incompressible lattice, and the enthalpy parameter $\gamma_\A > -2,5$. Note that $\gamma_\A <  - 2,5$ entails a phase separation \cite{Landstorfer:2013ly}. The electrostatic potential in the solid phase is denoted as $\phi_\S(\vecx,t)$. \\

\paragraph{\textbf{Electro-neutrality}} 
\label{par:electro_neutrality}
The electro-neutrality condition of $\Omega_\A$, $\Omega_\E$ and $\Omega_\C$ can be obtained by an asymptotic expansion of the balance equations in the electrochemical double layer at the respective surface $\Sigma$. We only briefly recapture the central conclusions and refer to \cite{Newman:1965aa,Landstorfer2016187,Landstorfer01012017,C6CP04142F,C5CP03836G,Landstorfer_2019} for details on the modeling, validation and the asymptotic derivation. Most importantly, electro-neutrality implies (i) that the double layer is in thermodynamic equilibrium, (ii) that there exists a potential jump through the interface $\Sigma_\iAE$ and (iii) that for monovalent electrolytes the cation mole fraction (or number density) is equal to the anion mole fraction, \ie 
\begin{align}
y_\EC = y_\EA =: y_\E ~.
\end{align}
Hence total number density $n_\Etot= n_\ES + n_\EC + n_\EA$  electrolyte concentration $n_\EC$ in terms of $y_\E$ write with eq. \eqref{eq:incomp_constraint} as 
\begin{align}
    n_\Etot & = n_\ES \cdot \frac{1}{1+2(\kappa_\E-1) y_\E} = n_\Etot(y_\E)   \\
    n_\EC   &= 
    n_\ES \frac{y_\E}{1+2(\kappa_\E-1) y_\E} = n_\EC(y_\E).
\end{align}

\subsubsection{Transport equation relations} 
\label{sub:transport}
    In the electrolyte $\Omega_\E$ we have two balance equations determining the concentration $n_\EC(\vecx,t)$ (or mole fraction $y_\E(\vecx,t)$) and the electrostatic potential $\phi_\E(\vecx,t)$ in the electrolyte \cite{Newman:1973ab,Newman:1973aa,Fuoss1978,doi:10.1021/j150551a038,LIU2014447,doi:10.1002/bbpc.19650690712}, \ie 
    \begin{align}
        \label{eq:transport_E_1}
        \frac{\partial n_\EC}{\partial t} &= -\div \flux_\EC \quad \text{with} \quad & \flux_\EC &= -   D_\E \cdot n_\Etot\,  \Gamma_\E^\tf \cdot \nabla y_\E + \frac{t_\EC}{e_0} \flux_\Eq     \\
                \label{eq:transport_E_2}
        0 &= -\div \flux_\Eq    \quad \text{with} \quad & \flux_\Eq &=  - S_\E \cdot n_\Etot\, \Gamma^\tf \nabla y_\E - \Lambda_\E n_\E \nabla \phi_\E
    \end{align}
    with (dimensionless) thermodynamic factor \cite{C3CP44390F}
    \begin{align}
        \Gamma_\E^\tf 
                                  &= y_\E \frac{\partial f_\E}{\partial y_\E} = 1 + 2 \kappa_\E \frac{y_\E}{1 - 2 y_\E} = \Gamma_\E^\tf(y_\E).
    \end{align}
%
    For simple Nernst--Planck-flux relation for the cation and anion fluxes \cite{C3CP44390F,sanfeld1968introduction}, \ie 
    \begin{align}
        \flux_\alpha = D_\alpha^\text{NP} \frac{n_\alpha}{\kT} ( \nabla \mu_\alpha - \frac{m_\alpha}{m_0} \nabla \mu_\ES  + e_0 z_\alpha n_\alpha \nabla \phi_\E) \qquad \alpha = \EA,\EC~, 
    \end{align}
    with constant diffusion coefficients $D_\EA^{\text{NP}}$ for the anion and $D_\EC^{\text{NP}}$ for the cation, we obtain (in the electro-neutral electrolyte)
    \begin{align}
    D_\E &= \frac{2 D_\EC^{\text{NP}} \cdot D_\EA^{\text{NP}}}{D_\EA^{\text{NP}} + D_\EC^{\text{NP}}}    & \quad        t_\EC &= \frac{D_\EC^{\text{NP}}}{D_\EA^{\text{NP}} + D_\EC^{\text{NP}}}                \\
    \Lambda_\E &= \frac{e_0^2}{\kT}(D_\EA^{\text{NP}} + D_\EC^{\text{NP}})  & S &= e_0(D_\EC^{\text{NP}} - D_\EA^{\text{NP}})   ~.
    \end{align}
However, for general Maxwell-Stefan type diffusion \cite{doi:10.1021/j150551a038,Fuoss1978,doi:10.1002/bbpc.19650690712,LIU2014447,Kim01012016} or cross-diffusion coefficients\cite{Latz:2011aa,dGM84} in the cation and anion fluxes, more \emph{complex} representations of the transport parameters $(t_\EC, S_\E, D_\E,\Lambda_\E)$ are obtained. In general, three of the transport parameters are \emph{independent}, and $S_\E$, $t_\EC$ and $\Lambda_\E$ are connected via
    \begin{align}
        \label{eq:relation_transport_parameter}
        \frac{\kT}{e_0} (2 t_C - 1) = \frac{S_\E}{\Lambda_\E} ~.
    \end{align}
    Further, $(t_\EC, S_\E, D_\E,\Lambda_\E)$ depend in general non-linearly on the electrolyte concentration $n_\EC$. However, it is sufficient for the sake of this work to assume constant values for the transport parameters $(t_\EC, S_\E, D_\E,\Lambda_\E)$, together with relation \eqref{eq:relation_transport_parameter}. \\
    
In the active particle $\Omega_\A$ we have two balance equations determining the concentration $n_{\A_C}(\vecx,t)$ (or mole fraction $y_\A$) and the electrostatic potential $\phi_\S(\vecx,t)$ in the solid phase $\Omega_\S$, \ie
    \begin{align}
        \label{eq:diffusion_active}
        \frac{\partial n_{\A_C}}{\partial t} &= -\div \flux_{\A_C} \quad \text{with} \quad & \flux_{\A_C} &= -  D_\A \cdot n_\Aell\, \Gamma_\A^\tf \cdot \nabla y_\A\\
        \label{eq:charge_balance_active}
        0 &= -\div \flux_\Sq    \quad \text{with} \quad & \flux_\Aq &=  - \sigma_\S(\vecx) \nabla \phi_\S
    \end{align}
    and (dimensionless) thermodynamic factor
    \begin{align}
        \Gamma_\A^\tf =  y_\A \frac{\partial f_\A}{\partial y_\A} = 1 + \frac{y_\A}{1- y_\A} + 2 \gamma_\A y_\A = \Gamma_\A^\tf(y_\A) ~.
    \end{align}
    The (solid state) diffusion coefficient $D_\A$ is further assumed to be 
    \begin{align}
        D_\A = D_{\A,\rref} \cdot (1-y_\A)~, \qquad D_{\A,\rref} = \const, 
    \end{align}
    where the term $(1-y_\A)$ accounts for a reduced (solid state) diffusivity due to crowding\cite{Landstorfer_2019}. Note that the electrical conductivity $\sigma_\S$ is different in the active phase $\Omega_\A$ and the conductive phase $\Omega_\C$, which form $\Omega_\S = \Omega_\A \cup \Omega_\C$. We account for this as 
    \begin{align}
\sigma_\S(\vecx) = \begin{cases}
    \sigma_{\A}, &\text{ if }\vecx \in \Omega_\A\\
    \sigma_{\C}, &\text{ if }\vecx \in \Omega_\C
\end{cases}~,
    \end{align}
In principle $\sigma_\A$ can be dependent on the amount of intercalated ions, \ie $\sigma_\A = \sigma_\A(y_\A)$, but for the sake of this work we assume $\sigma_\A  = \const$ and $\sigma_\C = \const$.

    
The charge flux of electrons in the solid phase $\Omega_\S$ is described via 
\begin{align}
                0 &= -\div \flux_\Sq    \quad \text{with} \quad & \flux_\Sq &=  - \sigma_\S(\vecx) \nabla \phi_\S~.
\end{align}

\subsubsection{Reaction Rate and affinity} 
\label{par:variables}
 
At the the interface $\Sigma_\iAE$ the intercalation reaction 
\begin{align}
        \label{eq:reaction_intercalation}
        \ceE{Li+} + \ceA{e-} \cearrow \ceA{Li} + \kappa_\E \cdot \ceE{S}~
\end{align}
occurs, which is modeled on the basis of (surface) non-equilibrium thermodynamics\cite{Landstorfer_2019}. Hence, the (surface) reaction rate $\us R$ can in general be written with $\alpha \in [0,1]$ as  \cite{Landstorfer01012017,SANFELD1989L521,DEFAY1977498,C5CP03836G,dgm2018} 
\begin{align}
        \label{eq:reaction_rate_surface}
        \us R & = \us L \cdot g\left( \frac{1}{\kT} \us \lambda\right)   &   &\text{with} & g(z) &:= \left(\exp{\alpha \cdot z}         -  \exp{-(1-\alpha) \cdot z} \right)~,
\end{align}
where 
\begin{align}
    \label{eq:surf_affinity_1}
\us \lambda_\iAE = e_0 \hat \phi_\S|_\iAE - e_0 \hat \phi_\E|_\iAE + \kT \cdot ( f_\A^j(y_\A|_\iAE) - f_\E(y_\E|_\iAE) )
\end{align}
denotes the surface affinity of the reaction \eqref{eq:reaction_intercalation}, which is  \emph{pulled back}  through the double layer to the respective \emph{points} (in an asymptotic sense) outside of the double layer. The quantity $\hat \phi_\E := \phi_\E - E_{\ce{Li+},\E}$ 
denotes the electrolyte potential vs. metallic Li and $\hat \phi_\S^j := \phi_\S^j -  E_{\ce{Li+},\A}^j$ 
the electrode $^j$ potential vs. metallic Li\cite{Landstorfer_2019}. Note again that $y_\A|_\iAE$ denotes the evaluation of $y_\A$ at the interface $\Sigma_{\A,\E}$ and that the surface affinity \eqref{eq:surf_affinity_1} is dependent on the chemical potential (or the mole fraction) evaluated at the interface.

Note that the non-negative function $\us L$ in \eqref{eq:reaction_rate_surface} ensures a non-negative entropy production $\us r_{\sigma,R}$ due to reactions on the surface, \ie $\us r_{\sigma,R} = \us \lambda \cdot \us R  > 0$. For the sake of this work we assume $\us L = \const$ and refer to \cite{Landstorfer_2019} for a detailed discussion on concentration dependency. 

\subsubsection{Balance equations} 
\label{par:balance_equations}
Applying the above modeling procedure for $j = \An,\Cat$ yields the following equation system, 
\begin{align}
    \label{ML:bal_dim_1}
    \frac{\partial n_{\A_C}^j}{\partial t} & = -\div \flux_{\A_C}^j   & \text{with} \quad & 
    \flux_{\A_C}^j = -  D_\A^j \cdot n_{\A,\lat}^j \, \Gamma_\A^j \cdot \nabla y_\A^j && \vecx \in \Omega_\A^j~,  \\
    \label{ML:bal_dim_2}
     0 &= -\div \flux_\Sq^j         & \text{with} \quad & 
     \flux_\Sq^j =  - \sigma_\S^j(\vecx) \nabla \phi_\S^j && \vecx \in \Omega_\S^j~, \\
    \label{ML:bal_dim_3}
    \frac{\partial n_\EC}{\partial t} &= -\div \flux_\EC  & \text{with} \quad & 
     \flux_\EC = -   D_\E \cdot n_\Etot \,  \Gamma_\E \cdot \nabla y_\E + \frac{t_\EC}{e_0} \flux_\Eq && \vecx \in \Omega_\E~,\\
    \label{ML:bal_dim_4}
     0 &= -\div \flux_\Eq & \text{with} \quad &
     \flux_\Eq =  - S_\E \cdot n_\Etot\, \Gamma_\E \nabla y_\E - \Lambda_\E n_\EC \nabla \phi_\E     && \vecx \in \Omega_\E
\end{align}
where \eqref{ML:bal_dim_1} is the balance of intercalated lithium within the active phase, \eqref{ML:bal_dim_2} the charge balance in the solid phase and \eqref{ML:bal_dim_2}$_2$ the electron flux, \eqref{ML:bal_dim_3} the balance of lithium ions in the electrolyte phase and  \eqref{ML:bal_dim_4} the charge balance in the electrolyte, where \eqref{ML:bal_dim_4}$_2$ is the flux of the charge $q_\E$ in the electrolyte. Note that $\sigma_\S^j(\vecx)$ incorporates the fact that the conductivity is far larger in the conductive additive phase $\Omega_\C^j$ then in the active particle phase $\Omega_\A^j$. The index $^j$ is necessary because anode and cathode are in general different materials, hence having different parameters and material functions, but \eqref{ML:bal_dim_1} yields a compact typeface.

The boundary conditions at the interface $\Sigma_\iAE^j = \Omega^j_\A \cap \Omega^j_\E$, where the intercalation reaction \ce{Li+ + e- <=> Li} occurs, read
\begin{align}    
    \label{ML:bal_dim_5}
    \flux_{\A_C}^j \cdot \vecn      & =     -\us L^j \cdot  g\big(\frac{1}{\kT} \us \lambda_\iAE^j \big)            && \text{on } \Sigma_\iAE^j   \\
    \flux_\Sq^j \cdot \vecn     &= - e_0 \us L^j \cdot  g\big(\frac{1}{\kT} \us \lambda_\iAE^j \big)            && \text{on } \Sigma_\iAE^j   \\
    \flux_\EC \cdot \vecn        & =      \us L^j\cdot   g\big(\frac{1}{\kT} \us \lambda_\iAE^j \big)            && \text{on } \Sigma_\iAE^j ~,  \\
        \label{ML:bal_dim_8}
    \flux_\Eq \cdot \vecn     &=  e_0 \us L^j \cdot  g\big(\frac{1}{\kT} \us \lambda_\iAE^j \big)            && \text{on } \Sigma_\iAE^j  
\end{align}
where by convention $\vecn$ is pointing from $\Omega_\A^j$ into $\Omega_\E^j$. At the interface $\Sigma_\iCE^j = \Omega^j_\C \cap \Omega^j_\E$ we have homogenous Neumann boundary conditions, \ie
\begin{align}    
            \label{ML:bal_dim_9}
    \flux_{\A_C}^j \cdot \vecn      & = 
    \flux_\Sq^j \cdot \vecn     =  
    \flux_\EC\cdot \vecn = \flux_\Eq \cdot \vecn         = 0   && \text{on } \Sigma_\iCE^j~.
\end{align}
Further we have \emph{global} boundary conditions for the electric current density $i$ leaving the anode (discharge, $i>0$) or entering the anode (charge, $i<0$), \ie
\begin{align}
                \label{ML:bal_dim_10}
    \flux_\Sq^\Cat \cdot \vecn &= - i && \text{on } \Sigma_\iCD^\Cat
\end{align}
as well as a reference value for the electrostatic potential, which reads 
\begin{align}
                \label{ML:bal_dim_11}
    \phi_S^\An &= 0 && \text{on } \Sigma_\iCD^\An~.
\end{align}

\subsection{Homogenization}
\label{ssub:homogen}
\subsubsection{Introduction} 
\label{ssub:introduction}
\begin{figure}[h] 
   
    \centering
    \includegraphics[width=\textwidth]{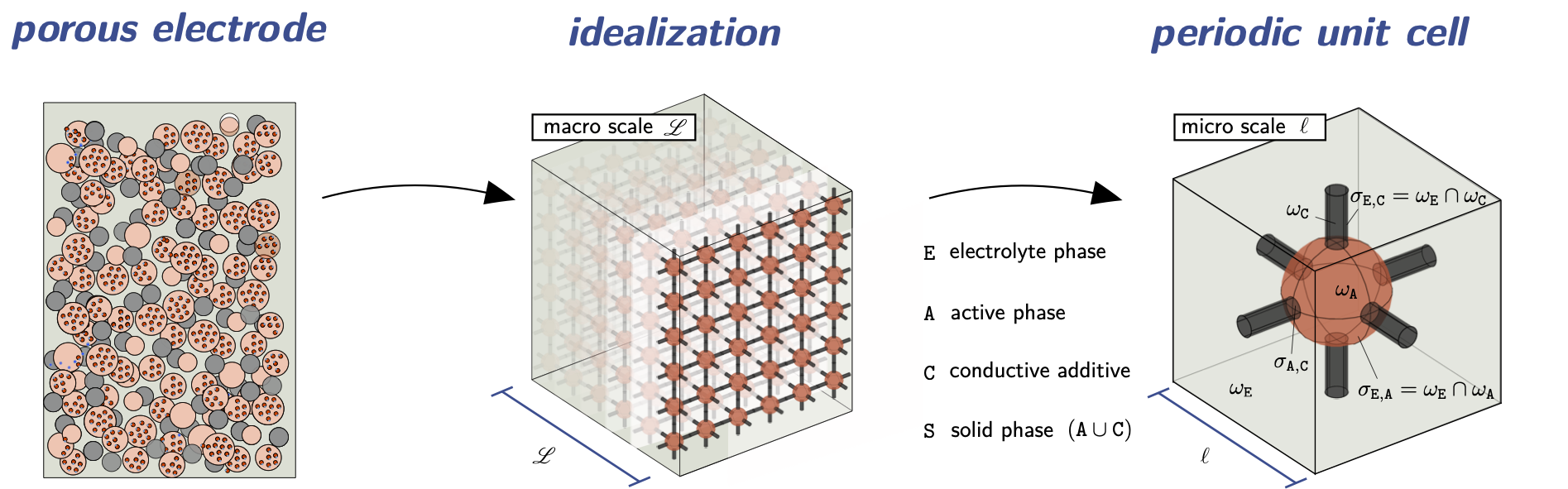}
    \caption{Sketch of the homogenization procedure. The porous electrode is simplified as a network of interconnected active phase spheres, yielding a unit cell $\omega$ containing one electrode particle.  \label{fig:electrode_hom}}
\end{figure}

An important feature of the coupled transport equation system \eqref{ML:bal_dim_1} - \eqref{ML:bal_dim_4} is the circumstance that the solid state diffusion $D_\A^j$ is far smaller than the electrolyte diffusivity $D_\E$. This accompanies, however, with smaller diffusion pathways on the length scale $\ell$ of the intercalation particles. The diffusivity of $Li$ in \ce{LiCoO2} is, for instance, about $D_\A \approx 10^{-12} \units{\frac{\Ucm}{\Us}}$, while the diffusivity of \ce{Li+} in DMC is in the order of $D_\E \approx \approx 10^{-5}\units{\frac{\Ucm}{\Us}}$. The macro-length scale is $\Ell \approx 1 \units{\Umum}$ while the micro-scale is $\ell \approx 1 \units{\Unm}$ (see Fig. \ref{fig:electrode_hom}). Hence the length scale ration $\frac{\ell}{\Ell} = \eps \approx 10^{-3}$ and  $D_\A^j \approx \eps^2  D_\E$. This motivates the re-scaling 
\begin{align}
    D_\A^j = \eps^2 \cdot D_\A^{j,\eps} && j = \An,\Cat,
\end{align}
which is subsequently exploited in the homogenization procedure. 

\subsubsection{Non-dimensionalization} 
\label{eq:subsection_non-dim_1}
For the sake of the homogenization as well as parameter studies and numerical implementations, it is necessary to non-dimensionalize the overall equation system. This is performed in several steps, starting from an initial non-dimensionalization for the homogenization, briefly sketched here. Consider the dimensional equation system
\begin{align}
    \frac{\partial w}{\partial t} &= \div D \nabla w          && \vecx \in \Omega_\E    \\
    D \nabla w \cdot \vecn &= R                              && \vecx \in \Sigma_\iAE
\end{align}
For the sake of the homogenization, it is convenient to introduce the scaling 
\begin{align}
    w &= u \cdot n              &   \tilde D & =  \frac{T}{\Ell^2} \cdot D                    \\
    t &= \tau \cdot T             &   \tilde R &=  \frac{T}{n} \frac{1}{\ell} \cdot R         \\
    x &= \xi \cdot \Ell           &   \eps    &= \frac{\ell}{\Ell}
\end{align}
which yields 
\begin{align}
    \frac{\partial u}{\partial \tau} &= \div_\xi D \nabla_\xi u          && \xi \in \tilde \Omega_\E \\
    \tilde D \nabla_\xi u \cdot \vecn &= \eps \tilde R                              && \xi \in \tilde \Sigma_\iAE.
\end{align}

\subsubsection{General Homogenization Framework} 
\label{ssub:general_homogenization_framework}
For a single electrode, dropping the super-index $^j$ and denoting in this subsection $x$  as non-dimensional space variable, we have essentially a coupled equation system of the form\footnote{Actually in the electrolyte phase we have two coupled equations, but for simplicity we sketch in the derivation here only a single (non-linear) equation.}
\begin{align}
    \label{eq:hom_intro_1}
        \frac{\partial u_\E}{\partial \tau}  &= \div \big(\tilde D_\E h_{\E}(u_\E) \nabla u_\E\big)                     && \in \Omega_\E        \\
    \frac{\partial u_\A}{\partial \tau}  &= \div \big(\eps^2 \tilde D_\A h_{\A}(u_\A) \nabla u_\A \big)                 && \in \Omega_\A                \\      
    \frac{\partial u_\S}{\partial \tau}  &= \div \big(\tilde D_\S^\eps h_{\S}(u_\S) \nabla u_\S \big)                   && \in \Omega_\S                \\      
                h_{\E}(u_\E) \tilde D_\E \nabla u_\E  \cdot \vecn &= \eps \tilde R_\E (u_\E \big|_{\Sigma_\iAE}^+ ,u_\A \big|_{\Sigma_\iAE}^-,u_\S \big|_{\Sigma_\iAE}^-)       && \text{on } \Sigma_\iAE       \\
                h_{\A}(u_\A) \eps^2 \tilde D_\A \nabla u_\A  \cdot \vecn &= \eps \tilde R_\A (u_\E \big|_{\Sigma_\iAE}^+ ,u_\A \big|_{\Sigma_\iAE}^-,u_\S \big|_{\Sigma_\iAE}^-)        && \text{on } \Sigma_\iAE       \\
                h_{\S}(u_\A) \tilde D_\S \nabla u_\S  \cdot \vecn &= \eps \tilde R_\S (u_\E \big|_{\Sigma_\iAE}^+ ,u_\A \big|_{\Sigma_\iAE}^-,u_\S \big|_{\Sigma_\iAE}^-)       && \text{on } \Sigma_\iAE       \\      
    \label{eq:hom_intro_end}
                h_{\E}(u_\E) \tilde D_\E \nabla u_\E  \cdot \vecn &=
                h_{\A}(u_\A) \eps^2 \tilde D_\A \nabla u_\A  \cdot \vecn = 
                h_{\S}(u_\A) \tilde D_\S \nabla u_\S  \cdot \vecn = 0   && \text{on } \Sigma_\iCE               
\end{align}
where $u_i, i = \A,\E,\S$, denotes the respective phase variable, $(\tilde D_i^\eps, \tilde R_i^\eps)$ the $\eps$-dependent non-equilibrium parameter, and $h_{i,k}$ captures non-linearities. Note that we abbreviate 
\begin{align}
        u_i \big|_{\Sigma_\iAE}^+ =: u_\E|^+ \quad \text{and} \quad u_i \big|_{\Sigma_\iAE}^- = : u_i|^-~.
\end{align}
We consider a multi-scale expansion ($i=\A,\E,\S$)
\begin{align}
        u_i(\vecx,t) = \sum_{k=0}^{\infty} \eps^k u_i^k(x,y,t) \quad \text{with} \quad y = \frac{x}{\eps} 
\end{align}
whereby
        \begin{align}
                \nabla &= \nabla_x + \eps^{-1} \nabla_y\\
                \div &= \div_x + \eps^{-1} \div_y~.              
        \end{align}
For non-linear functions $h=h(u)$, we consider the $\eps-$Taylor expansion
\begin{align}
        h(u) = \sum_{k=0}^{\infty} \frac{1}{k!} \frac{d^k h}{d u^k}(u-u^0)  
                                                &= h(u^0) + \eps u^1 h'(u^0) + \O(\eps^2)
\end{align}
and for $g=g(u_\E^+,u_\A^-,u_\S^-)$ a multi-dimensional Taylor expansions, \ie
\begin{align}
        g(u_\E|^+,u_\A|^-,u_\S^-) = g( u_\E^0|^+,u_\A^0|^-,u_\S^0|^-) + \sum_{i=\E,\A,\S} \eps \cdot \partial_{u_i} g \big|_{u_\E^0|^+,u_\A^0|^-,u_\S^0|^-} \cdot u_i^1|^\pm + \O(\eps^2)~.
\end{align}
We abbreviate 
\begin{align}
        &h^0 :=  h(u^0)\, , ~   h^1 := u_1 \cdot h'(u^0) ~, \quad \text{and} \quad g^0 := g( u_\E^0|^+,u_\A^0|^-,u_\S^0|^-) \\
    &g^1 := u_\E^1 \partial_{u_\E} g (u_\E^0|^+,u_\A^0|^-,u_\S^0|^-) + u_\A^1 \partial_{u_\A} g (u_\E^0|^+,u_\A^0|^-,u_\S^0|^-)  + u_\S^1 \partial_{u_\S} g (u_\E^0|^+,u_\A^0|^-,u_\S^0|^-)  
\end{align}
and expand thus
\begin{align}
        h_i &= h_i^0 + \eps h_i^1                       \qquad\qquad    g = g^0 + \eps g^1~.
\end{align}
Insertion of the multi scale expansion yields essentially a sequence of PDEs in the orders of $\eps$.\\

\paragraph{\textbf{Order $\eps^{-2}$:}} 
\label{ssub:order_eps_2}
For $i=\E,\S$ we have 
\begin{align}
        \div_y \tilde D_i h_{i}^0  \nabla_y u_i^0 
\end{align}
with 
\begin{align}
         \tilde D_i h_{i}^0 \nabla_y u_i^0  = 0
\end{align}
This yields $u_i^0 = u_i^0(\vecx,t)$ by the maximum principle.\\

\paragraph{\textbf{Order $\eps^{-1}$:}} 
\label{ssub:order_eps_1}
For $i=\E,\S$ we have due to $u_i^0 = u_i^0(\vecx,t)$ we have 
\begin{align}
\div_y \tilde D_i h_{i}^0 \nabla_y  u_i^1        =  0
\end{align}
with
\begin{align}
        (\tilde D_i h_{i}^0 \nabla_x u_i^0      +  \tilde D_i h_{i}^0 \nabla_y u_i^1) \cdot \vecn = 0~.
\end{align}
This yields essentially 
\begin{align}
        u_i^1(x,y,t) = - \nabla_x u_i^0 \cdot \vecchi_i(y)
\end{align}
where $\vecchi_\E = (\chi_\E^1,\chi_\E^2,\chi_\E^3)$ satisfies the \textbf{cell problem } ($k=1,2,3$)
\begin{align}
    \label{eq:CP_E}
(\textbf{CP}_\E) \qquad \begin{cases}
        \div_y \nabla_y \chi_\E^k & = 0 \qquad          \vecy \in \omega_\E, \, i =\E,\S   \\
        \nabla \chi_\E^k \cdot \vecn & = n_k    \qquad \text{on} \sigma_\iAE \\
        \quad \chi_\E^k & \text{ is periodic } 
\end{cases}
\end{align}
and $\vecchi_\S = (\chi_\S^1,\chi_\S^2,\chi_\S^3)$ satisfies the \textbf{cell problem } ($k=1,2,3$)
\begin{align}
        \label{eq:CP_S}
(\textbf{CP}_\S) \qquad \begin{cases}
        \div_y \left( h_{3,\S}(\vecy) \cdot \left( \vece^k - \nabla_y \chi_\S^k   \right) \right)       & = 0    \qquad  \vecy \in \omega_\S,  \\ 
         h_{3,\S}(\vecy)  \cdot \left( \vece^k - \nabla_y \chi_\S^k   \right) \cdot \vecn & = 0  \qquad \text{on} \sigma_\iAE  \\
        \quad \chi_\S^k & \text{ is periodic } 
\end{cases}
\end{align} \\
\paragraph{\textbf{Order $\eps^{0}$}} 
\label{ssub:order_eps_0}
Since $u_\E^0 = u_\E^0(\vecx,t)$ we have for $i=\E,\S$
\begin{align}
        \label{eq:Ntilde R_eps0_1} 
         \div_x  \tilde D_i h_{i}^0 \nabla_x u_i^0       +  \div_x  \tilde D_i h_{i}^0 \nabla_y u_i^1    + \div_y \tilde D_i h_{i}^0 \nabla_x u_i^1      + \div_y \tilde D_i h_{i}^1 \nabla_y  u_i^1      = \frac{\partial u_i^0}{\partial \tau}         
\end{align}
and
\begin{align}
         \div_y \tilde D_\A h_{\A}^0  \nabla_y u_\A^0 =    \frac{\partial u_\A^0}{\partial \tau} 
\end{align}
with 
\begin{align}
\Big(\tilde D_\E h_{\E}^0 \nabla_x u_\E^1       +  \tilde D_\E h_{\E}^1 \nabla_x u_\E^0         +  \tilde D_\E h_{\E}^1 \nabla_y u_\E^1         \Big) \cdot \vecn& =  \tilde R_\E^0\\
\Big(\tilde D_\S h_{\S}^0 \nabla_x u_\S^1       +  \tilde D_\S h_{\S}^1 \nabla_x u_\S^0         +  \tilde D_\S h_{\S}^1 \nabla_y u_\S^1         \Big) \cdot \vecn& =  \tilde R_\S^0\\
\Big(\tilde D_\A h_{\A}^0 \nabla_y u_\A^0  \Big) \cdot \vecn & = \tilde R_\A^0 
\end{align}
 
 
Eq. \eqref{eq:Ntilde R_eps0_1} leads (by an integration over $\omega_\E$ and integration by parts) for $i = \E$ to 
\begin{align}
 \frac{\partial u_\E^0}{\partial \tau}   
                 & =     \div_x \big( \tilde D_\E h_{\E}^0      \tenpi_\E \cdot  \nabla_x u_\E^0 \big)   +  \frac{a_\iAE}{\psi_\E }  \cdot \frac{1}{\area{\sigma_\iAE} }        \int_{\sigma_\iAE} \tilde R_\E^0(u_\E^0, u_\S^0,u_\A^0|_{\partial \omega_\A}) dA(\vecy) 
\end{align}
with
\begin{align}
        \tenpi_\E &:= \frac{1}{\vol{\omega_\E} } \int_{\omega_\E} (\tenId -  \tennabla_y \vecchi_\E ) dV(y)  =  \tenId -  \frac{1}{\vol{\omega_\E} } \int_{\omega_\E}\tennabla_y \vecchi_\E  dV(y)      \\      
        \psi_\E &= \frac{\int_{\omega_\E} 1 \,dV }{\int_{\omega} 1 \,dV  } = \frac{\vol{\omega_\E}}{\vol{\omega} }  \\  
        a_\iAE &= \frac{\int_{\sigma_\iAE} 1\,dA}{{\int_{\omega} 1 \,dV } }  = \frac{\area{\sigma_\iAE}}{\vol{\omega}}
\end{align}
Note that $u_\A^0=u_\A^0(\vecx,\vecy,t)$ and that $u_\A^0|_{\partial \omega_\A}$ denotes an evaluation of $u_\A^0$ at the boundary of the micro-domain $\omega_\A$. Eq. \eqref{eq:Ntilde R_eps0_1} leads (by an integration over $\omega_\S$) for $i = \S$ to 
\begin{align}
\frac{\partial u_\S^0}{\partial \tau}     
                 & =     \div_x \big( \tilde D_\S^0 h_{\S}^0    \tenpi_\S \cdot  \nabla_x u_\S^0 \big)   +  \frac{a_\iAE}{\psi_\S }  \frac{1}{\area{\sigma_\iAE} }      \int_{\sigma_\iAE} \tilde R_\S^0(u_\E^0, u_\S^0,u_\A^0|_{\partial \omega_\A}) dA(\vecy) 
\end{align}
with
\begin{align}
        \tenpi_\S &:= \frac{1}{\vol{\omega_\S} } \int_{\omega_\S} h_{\S,3} (y) (\tenId -  \tennabla_y \vecchi_\S ) dV(y)  \\    
        \psi_\S &= \frac{\int\limits_{\omega_\S} 1 \,dV }{\int\limits_{\omega} 1 \,dV  } =  \frac{\vol{\omega_\S}}{\vol{\omega} }  
\end{align}
%
Hence we obtain for the equation system     \eqref{eq:hom_intro_1} --  \eqref{eq:hom_intro_end} via periodic homogenization the coupled micro-macro balance equation system 
\begin{align}
\frac{\partial u_\E^0(\vecx,t)}{\partial \tau}    & =  \div_x \big( \tilde D_\E h_{\E}^0        \tenpi_\E \cdot  \nabla_x u_\E^0 \big)   +      \frac{a_\iAE}{\psi_\E}  \frac{1}{\area{\sigma_\iAE} }   \int_{\sigma_\iAE} \tilde R_\E^0(u_\E^0, u_\S^0,u_\A^0|_{\partial \omega_\A}) dA(\vecy) \\
\frac{\partial u_\S^0(\vecx,t)}{\partial \tau}    & =  \div_x \big( \tilde D_\S^0 h_{\S}^0      \tenpi_\S \cdot  \nabla_x u_\S^0 \big)   +      \frac{a_\iAE}{\psi_\S}  \frac1{\area{\sigma_\iAE} }     \int_{\sigma_\iAE} \tilde R_\S^0(u_\E^0, u_\S^0,u_\A|_{\partial \omega_\A}) dA(\vecy) \\
\frac{\partial u_\A^0(x,y,t)}{\partial \tau}      & =   \div_y \tilde D_\A h_{\A}^0  \nabla_y u_\A^0 
\end{align}
with 
\begin{align}
        \Big(\tilde D_\A h_{\A}^0 \nabla_y u_\A^0  \Big) \cdot \vecn & = \tilde R_\A^0~.
\end{align}\\

\paragraph{\textbf{Spherical particles}} 
\label{par:spherical_particles}
In the following, we assume on the \textbf{macro-scale} a 1-D approximation $x \in [0,1]$ as well as \textbf{spherical particles} $\omega_\A$, \ie $\omega = [-0.5,0.5]^3$ ($\vol{\omega} = 1$), $\omega_\A = B_{\tilde r_\A}(\mathbf{0})$, $ \tilde r_\A < 0.5$, $\omega_\E = \omega \backslash \omega_\A$, yielding  
\begin{align}
    \label{eq:_hom_pde_1}
\psi_\E \frac{\partial u_\E^0(x,t)}{\partial \tau}        & =  \partial_x \big( \psi_\E \tilde D_\E h_{\E}^0    \tenpi_\E \cdot   \partial_x u_\E^0 \big)\!+\!\theta_\iAE \tilde R_\E^0\left(u_\E^0,u_\S^0,u_\A^0(x,\tilde r_\A,t) \right) ~, &&  x \in [0,1]           \\
    \label{eq:_hom_pde_2}
\psi_\S \frac{\partial u_\S^0(x,t)}{\partial \tau}        & =  \partial_x \big( \psi_\S \tilde D_\S^0 h_{\S}^0  \tenpi_\S \cdot   \partial_x u_\S^0 \big)\!+\!\theta_\iAE \tilde R_\S^0\left(u_\E^0,u_\S^0,u_\A^0(x,\tilde r_\A,t) \right)      ~, && x \in [0,1]               \\
     \label{eq:_hom_pde_3}
\frac{\partial u_\A^0(x,r,t)}{\partial \tau}  & =       \frac{1}{r^2} \partial_r \big( r^2 \tilde D_\A h_{\A}^0  \partial_r u_\A^0 \big)        ~, && \hspace{-1.8cm}    x \in [0,1]~, r \in (0,\tilde r_\A)    
\end{align}
with 
\begin{align}
                \Big( \tilde D_\A h_{\A}^0 \partial_r u_\A^0 \Big) \Big|_{r=r_\A} & = \tilde R_\A^0 \left(u_\E^0(x,t),u_\S^0(x,t),u_\A^0(x,\tilde r_\A,t)\right) 
\end{align}
and
\begin{align}
        \theta_\iAE = \frac{\area{\sigma_\iAE}}{\vol{\omega}} = 4 \pi \tilde r_\A^2 - \area{\sigma_{\A,\C}}~.
\end{align}
\begin{figure}[h] 
    \centering
    \includegraphics[width=\textwidth]{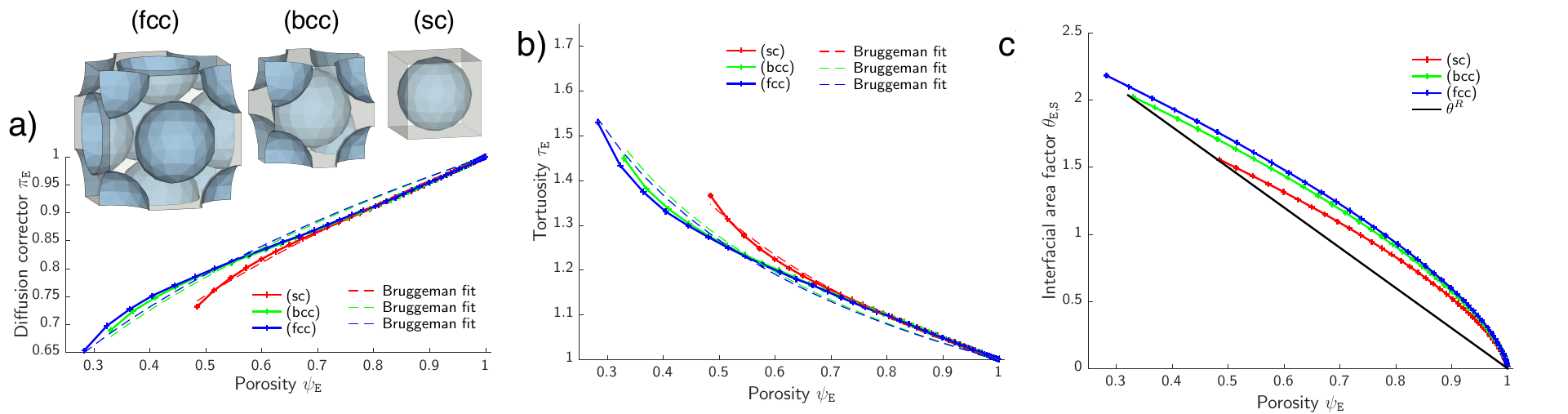}
    \caption{Porous media parameters for simple cubic, body centered cubic and face centered cubic micro-structures (Fig. 15 of \cite{LANDSTORFER2021110071}, reprinted with permission of Elsevier).}
        \label{fig:SC}
\end{figure}

For spherical particles, various possibilities regarding their micro-structural packing arise, most prominent (i) simple cubic, (ii) body centered cubic and (iii) face centered cubic, see Fig. \ref{fig:SC}. For a given micro-structure (or periodic unit cell $\omega$), the porous media parameters $(\psi_\E,\tenpi_\E,\theta_\iAE)$ and $(\tenpi_\S,\psi_\S)$ can (numerically) be computed by solving the cell problems \eqref{eq:CP_E} and \eqref{eq:CP_S}, respectively. Note that for equally sized spherical particles the actual micro-structure is exclusively encoded in the porous media parameters $(\psi_\E,\tenpi_\E,\theta_\iAE,\tenpi_\S,\psi_\S)$ of the homogenized equation system \eqref{eq:_hom_pde_1}-\eqref{eq:_hom_pde_3}. Note that, quite commonly, the (scalar) tortuosity corrector $\tau_\E$ is introduced via an  effective diffusion coefficient\cite{Chung_2013,Tjaden:2016wk,Newman:1973aa,Ebner_2014},  which simply yields in our notation $\tau_\E = (\pi_\E)^{-1}$ \cite{LANDSTORFER2021110071}.

We refer to \cite{LANDSTORFER2021110071} for more complex micro-structure geometries as well as their meshing and numerical solution of \eqref{eq:CP_E}, and proceed for the sake of this work with a simple cubic micro-structure.

\subsubsection{Homogenized Battery Model} 
\label{ssub:homogenized_battery_model}
Applying this homogenization scheme to the equation system \eqref{ML:bal_dim_1} - \eqref{ML:bal_dim_4} and its boundary conditions  \eqref{ML:bal_dim_5} --  \eqref{ML:bal_dim_9}, dropping the leading order index $^0$, and reinserting the scalings of section \ref{ssub:introduction} and \ref{eq:subsection_non-dim_1} yields the following equation system ($j = \An,\Cat$):
\begin{align}  
    \frac{\partial n_{\A_C}^j}{\partial t} & = - \frac{1}{r^2} \partial_r J_{\A_C}^j   && (r,x) \in [0,r_\A^j] \times I^j \nonumber \\
 \label{ML:bal_dim_1_hom}
    & \text{with} \quad  
    J_{\A_C}^j = -  D_\A^j \cdot n_{\A,\lat}^j \, \Gamma_\A^j \cdot r^2 \,  \partial_r y_\A^j \\
     0 &= -\partial_x J_\Sq^j   - e_0 \theta_\iAE \frac{1}{\ell} \us{L}^j \cdot   g     && x \in I^j~,  \nonumber  \\
     \label{ML:bal_dim_2_hom}
     & \text{with} \quad  
     J_\Sq^j =  - \psi_\S \tenpi_{\S}^\sigma \sigma_\C^j \partial_x \phi_\S^j \\ 
    \psi_\E \frac{\partial n_\EC}{\partial t} &= -\partial_x J_\E  + \theta_\iAE \frac{1}{\ell} \us{L}_\E \cdot   g && x \in I~,\nonumber \\
      \label{ML:bal_dim_3_hom}
    & \text{with} \quad    J_\EC = -  \psi_\E \tenpi_\E D_\E \cdot n_\Etot \,  \Gamma_\E \cdot \partial_x y_\E + \frac{t_\EC}{e_0} J_\Eq \\
     0 &= -\partial_x J_\Eq +  e_0  \theta_\iAE \frac{1}{\ell} \us{L}^j\cdot   g && x \in I \nonumber \\
         \label{ML:bal_dim_4_hom}
     & \text{with} \quad  J_\Eq =  - \psi_\E \tenpi_\E S_\E \cdot n_\Etot\, \Gamma_\E \nabla y_\E - \psi_\E \tenpi_\E \Lambda_\E n_\EC \nabla \phi_\E     
\end{align}
where
\begin{align}
    I^\An &= [0,W^\An], ~ I^\Sep = [W^\An,W^\An\!+\!W^\Sep], ~ I^\Cat = [W^\An\!+\!W^\Sep,W] \\
    W &= W^\An + W^\Sep + W^\Cat,~ I = I^\An \cup I^\Sep \cup I^\Cat = [0,W], \\
    g &= g\big(y_\E(x,t),\phi_\E(x,t),\phi_S(x,t),y_\A(x,r_\A^j,t)\big)~,
\end{align}
and boundary conditions 
\begin{align}
    J_\A^j \Big|_{r=r_\A} \!\!= - r_\A^2 \us{L}^j \cdot g~, \quad J_\A^j \Big|_{r=0}  \!\!=  0~, \quad  J_\Sq^\Cat\Big|_{x = W}  \!\!= -  i~, \quad \phi_S^\An \Big|_{x = 0}  = 0 ~.
\end{align}
Note that $t_\EC = \const$ allows us to rewrite \eqref{ML:bal_dim_3_hom} with \eqref{ML:bal_dim_4_hom} as 
\begin{align}
\psi_\E \frac{\partial n_\EC}{\partial t} &= -\partial_x J_\E  + \theta_\iAE \frac{1}{\ell} (1-t_\EC) \us{L}_\E \cdot   g && x \in I~,\nonumber \\          
    \label{ML:bal_non_dim_3b}
          & \text{with} \quad  J_\EC = -  \psi_\E \tenpi_\E D_\E \cdot n_\Etot \,  \Gamma_\E \cdot \partial_x y_\E                 ~,
\end{align}
which is further used.

\paragraph{\textbf{Cell Voltage, Capacity, C-Rate}} 
\label{par:cell_voltage_capacity_c_rate}
To introduce proper scalings an non-dimensionalizations, some definitions of important (global) quantities are required. For a lithium ion battery, the following quantities are of most importance: 
\begin{enumerate}
    \item the cell voltage 
    \begin{flalign}
     E &:= \phi_\S^\An|_{x=0} - \phi_\S^\Cat|_{x=W}              ~ && \uunits{\UV}~,
    \end{flalign}
    \item the basic (electrode) capacity ($j=\An,\Cat$)    
    \begin{flalign}
                Q^{i,0} &:= \int_{I^j}  \frac{4 \pi}{\ell^3} \int_0^{r_\A^j} e_0 n_{\A,\lat}^j r^2 \, dr  dx  
                 = W^i \psi_\A^j q_\A^j   = \const ~  && \uunits{\UAh \, \Um{^2}}~,
    \end{flalign}    
    \item the present (electrode) capacity ($j=\An,\Cat$)
    \begin{flalign}
        Q^i(t) &:= \int_{I^j}  \frac{4 \pi}{\ell^3} \int_0^{r_\A^j} e_0 n_{\A_C}^j r^2 \, dr  dx ~ &&\uunits{\UAh \,  \Um{^2}}~,
    \end{flalign}
    \item the (electrode) status of charge 
    \begin{flalign}
                \bar y_\A^j(t) & := \frac{Q^i(t)}{Q^{i,0}} 
                = \frac{1}{W^j} \int_{I^j}   \frac{3}{r_\A^3} \int_0^{r_\A^j}  y_\A^j \cdot r^2 \, dr  dx ~ && \uunits{1}~, 
    \end{flalign}
    \item and the $1$-C current density 
    \begin{flalign}
        i_C & :=   \frac{Q^{\Cat,0}}{1 \uunits{\Uh}}         ~ && \uunits{\UA \Um{^2}}~,
    \end{flalign}
    which yields for the current density the scaling 
    \begin{align}
        i = C_h \cdot i_C~,
    \end{align}
    where $C_h \in \mathR$ is the C-Rate. 
\end{enumerate}
Note that 
\begin{align}
            Q^\Cat(t)\!-\!Q^{\Cat,0}  & = \int_{I^\Cat} \frac{e_0  4 \pi}{(\ell^\Cat)^3}  \int_0^{r_\A^\Cat}  \int_0^t r^2 \partial_ t n_{\A_C}^\Cat  \, dt \, dr \, dx 
                      = \int_0^t i  \, dt ~.
\end{align}
For a constant current discharge ($i = \const$) we obtain thus 
\begin{align}
                    \bar y_\A^\Cat(t) - \bar y_\A^\Cat(t=0) = \frac{C_h}{1 \uunits{\Uh}}  \cdot t~,
\end{align}
which introduces the time re-scaling $\tau = \frac{C_h}{1 \uunits{\Uh}}  \cdot t$. 

\paragraph{\textbf{Scalings}} 
\label{par:scalings}
Summarized, we use the scalings 
\begin{align}
      t &=  \frac{1 \uunits{\Uh}}{C_h} \tau \quad,\; \tau \in [0,1]     %
&   i_C & =  \frac{Q^{\Cat,0}}{1 \uunits{\Uh}}       = \frac{W^\Cat \psi_\A^\Cat q_\A^\Cat}{1 \uunits{\Uh}}\\
   n_{\E,\tot} &= \tilde n_{\E,\tot} n_{\E,\rref}      %
&    x &= \xi W  \quad,\;\xi \in [0,1]                  \\
n_{\A_C}^j &= n_{\A,\lat}^j y_\A^j ~,\; y_\A^j\in[0,1]   %
& \phi &=  \frac{\kT}{e_0} \tilde \phi \\
    r &= r_\A \nu \quad,\; \nu \in [0,1]                      %
&  c_\E &= \frac{\partial n_\EC(y_\E)}{\partial y_\E} \frac{1}{n_{\E,\rref}} \\
n_\EC &= 
            \tilde n_{\E,\tot} \, n_{\E,\rref} \, y_\E  \quad,\; y_\E \in [0,1]   %
&      \eta_n^{\E,j} &:=  \frac{n_{\A,\lat}^j}{n_{\E,\rref}^j} \\
\eta_x^j & = \frac{W^j}{W} %
& & \\
\tilde r_\A^j &:= \frac{\ell^j}{r_\A^j} & \eta_W^\Cat &:= \frac{\psi_\A^\Cat \cdot W^\Cat}{W}  \\
%
        \tilde \sigma_\S^j       &:= \frac{1 \uunits{\Uh}}{W^2} \frac{\frac{\kT}{(e_0)^2}}{n_{\A,\lat}^j } \cdot \sigma_\S^j
    %
& \tilde \Lambda_\E             &=: \frac{1 \uunits{\Uh}}{W^2} \frac{\kT}{e_0^2} \Lambda_\E \\
\tilde D_\E     &:= \frac{1\uunits{\Uh} }{W^2} D_\E     
    & \tilde S_\E &:=  (2 t_C - 1) \tilde \Lambda_\E \\   
\ustilde{L}^j         &:= \frac{1}{n_{\A,\lat}^j} \big( {1 \uunits{\Uh}} \big)  \frac{1}{\ell^j} \us{L}^j %
&        \tilde D_{\A,\rref}^j & := \frac{1 \uunits{\Uh}}{(\ell^j)^2} \cdot D_{\A,\rref}^j 
\end{align}
and definitions 
\begin{align}
        \Omega^\Cat  &= [0,\eta_x^\Cat],        \Omega^\Sep = [\eta_x^\Cat,\eta_x^\Cat\!+\!\eta_x^\Sep]~, ~\Omega^\An =  [\eta_x^\Cat\!+\!\eta_x^\Sep,1], \Omega^\El = \Omega^\Cat \cup \Omega^\An \\ 
        \Omega_\A^\Cat & = [0,\eta_x^\Cat] \times [0,1]~, \Omega_\A^\An = [\eta_x^\Cat+\eta_x^\Sep,1] \times [0,1] \quad \Omega_\A = \Omega_\A^\Cat \cup \Omega_\A^\An
\end{align}
as well as

\begin{align}
        y_\A &:= \begin{cases}
                        y_\A^\Cat               &        (\xi,r) \in \Omega_\A^\Cat \\
                        y_\A^\An                &            (\xi,r) \in \Omega_\A^\An \\               
        \end{cases}
        &&& \vspace{-0.5cm}
        \hat \sigma_\S  &:= \begin{cases}
        \psi_\S^\Cat \tenpi_{\S,\sigma}^\Cat \tilde \sigma_\C^\Cat              &        \xi \in \Omega^\Cat \\
        \psi_\S^\An \tenpi_{\S,\sigma}^\An   \tilde \sigma_\C^\An                       &        \xi \in \Omega^\An \\
        \end{cases}\\
        \phi_\S &:= \begin{cases}
                        \phi_\S^\Cat    &        \xi \in \Omega^\Cat \\
                        \phi_\S^\An     &                \xi \in \Omega^\An \\          
        \end{cases}
        &&&\vspace{-0.5cm}
        \hat D_\A  &:= \begin{cases}
                \tilde D_{\A,\rref}^\Cat \, (1-y_\A) \Gamma_\A^\Cat(y_\A) \cdot \nu^2   &        (\xi,r) \in \Omega_\A^\Cat \\
                \tilde D_{\A,\rref}^\An \,  (1-y_\A) \Gamma_\A^\An(y_\A) \cdot \nu^2         &   (\xi,r) \in \Omega_\A^\An \\
        \end{cases}\\
        \psi_\E &:= \begin{cases}
                                        \psi_\E^\Cat    &        \xi \in \Omega^\Cat \\
                    \psi_\E^\Sep    &    \xi \in \Omega^\Sep \\
                                        \psi_\E^\An      &       \xi \in \Omega^\An \\
                                \end{cases} &&&\vspace{-0.5cm}
                                \hat D_\E  &:= \begin{cases}
                         \psi_\E^\Cat \tenpi_\E^\Cat \tilde D_\E \cdot \tilde n_{\E,\tot}(y_\E)  \,  \Gamma_\E(y_\E)  &  \xi \in \Omega^\Cat \\
                         \psi_\E^\Sep \tenpi_\E^\Sep \tilde D_\E \cdot \tilde n_{\E,\tot}(y_\E)  \,  \Gamma_\E(y_\E)  &  \xi \in \Omega^\Sep \\
                         \psi_\E^\An  \tenpi_\E^\An  \tilde D_\E \cdot \tilde n_{\E,\tot}(y_\E)  \,  \Gamma_\E(y_\E)  &  \xi \in \Omega^\An \\
                \end{cases}                 
                \\
                \psi_\S &:= \begin{cases}
                                        \psi_\S^\Cat &   \xi \in \Omega^\Cat \\
                                        \psi_\S^\An      &       \xi \in \Omega^\An \\
                                \end{cases}     &&&\vspace{-0.5cm}
        \hat S_\E  &:= \begin{cases}
                \psi_\E^\Cat \tenpi_\E^\Cat   \cdot \tilde S_\E \cdot \tilde n_{\E,\tot}(y_\E) \, \Gamma_\E(y_\E) &      \xi \in \Omega^\Cat \\
                \psi_\E^\Sep \tenpi_\E^\Sep   \cdot \tilde S_\E \cdot \tilde n_{\E,\tot}(y_\E) \, \Gamma_\E(y_\E) &      \xi \in \Omega^\Sep \\
                \psi_\E^\An  \tenpi_\E^\An    \cdot \tilde S_\E \cdot \tilde n_{\E,\tot}(y_\E) \, \Gamma_\E(y_\E) &      \xi \in \Omega^\An \\
        \end{cases} \\    
    \eta_n^\E &:= \begin{cases}
                        \eta_n^{\E,\Cat}    &    \xi \in \Omega^\Cat \\
                                0                                                                &       \xi \in \Omega^\Sep \\
                        \eta_n^{\E,\An}          &       \xi \in \Omega^\An \\
                \end{cases} &&&\vspace{-0.5cm}
        \hat \sigma_\E  &:= \begin{cases}
                \psi_\E^\Cat \tenpi_\E^\Cat  \tilde \Lambda_\E  \cdot \tilde n_{\E,\tot}(y_\E) \, y_\E &        \xi \in \Omega^\Cat \\
                \psi_\E^\Sep \tenpi_\E^\Sep  \tilde \Lambda_\E  \cdot \tilde n_{\E,\tot}(y_\E) \, y_\E &        \xi \in \Omega^\Sep \\
                \psi_\E^\An  \tenpi_\E^\An   \tilde \Lambda_\E  \cdot \tilde n_{\E,\tot}(y_\E) \, y_\E &        \xi \in \Omega^\An \\
        \end{cases}\\
    \theta &:= \begin{cases}
                        \theta_\iAE^\Cat &       \xi \in \Omega^\Cat \\
                                0                                                                &       \xi \in \Omega^\Sep \\
                        \theta_\iAE^\An          &       \xi \in \Omega^\An \\
                \end{cases}  &&&                        \vspace{-0.5cm}
                R
        &:= \begin{cases}
                         \ustilde L^\Cat g(\ustlambda^\Cat)      &       \xi \in \Omega^\Cat \\
                                0                                                                &       \xi \in \Omega^\Sep \\
                         \tilde L^\An g(\ustlambda^\An)          &       \xi \in \Omega^\An \\
                \end{cases}\\
        \tilde r_\A &:= \begin{cases}
                        \tilde r_\A^\Cat                &  \xi \in \Omega^\Cat \\
                        \tilde r_\A^\An                 &        \xi \in \Omega^\An \\          
        \end{cases}
\end{align}
which yields 
%
\begin{align}
    \label{ML:bal_dim_1_hom_nd}
    \nu^2 \tilde r_\A^2 C_h \frac{\partial y_\A}{\partial \tau} & = -  \partial_\nu \tilde J_{\A_C}   & \text{with} \quad & 
    \tilde J_{\A_C} \!=\! - \hat D_\A \, \partial_\nu y_\A  && (\nu,\xi) \in \Omega_\A \\
    \label{ML:bal_dim_2_hom_nd}    
     0 &= -\partial_\xi \tilde J_\Sq    -  \theta \cdot R     & \text{with} \quad & 
    \tilde J_\Sq \!=\! -  \hat \sigma_\S \partial_\xi \tilde \phi_\S && \xi \in \Omega^\El~, \\
    \label{ML:bal_non_dim_3_hom_nd}
     \psi_\E    C_h c_\E   \frac{\partial y_\E}{\partial \tau} &= -\partial_\xi \tilde J_\EC + \eta_n^\E (1\!-\!t_\EC) \theta \cdot R  & \text{with} \quad & 
          \tilde J_\EC \!=\! -  \hat D_\E(y_\E)  \partial_\xi y_\E && \xi \in [0,1]~, \\  
    \label{ML:bal_non_dim_4_hom_nd}    
     0 &= -\partial_\xi \tilde J_\Eq + \eta_n^\E \theta \cdot R & \text{with} \quad &
     \tilde J_\Eq \!=\!  - \hat S_\E(y_\E) \partial_\xi y_\E - \hat \sigma_\E \partial_\xi \tilde \phi_\E     && \xi \in [0,1]~.
\end{align}
The boundary conditions read 
%
\begin{align}
    \tilde J_{\A_C} \Big|_{\nu=1}  &= - \tilde r_\A R   ~, &  \tilde J_{\A_C} \Big|_{\nu=0}  &= 0 \\
    \hat \sigma_\C^\Cat \partial_\xi \tilde \phi_\S^\Cat \Big|_{\xi = 1}  &= C_h  \eta_W^\Cat ~,%
          &  \tilde \phi_S^\An \Big|_{\xi = 0} & = 0 ~,
\end{align}
with additional homogenous Neumann boundary conditions for all unspecified boundaries.\\

\paragraph{\textbf{Sign convention for the current density}} 
\label{par:sign_convention_for_the_current_density}
For the reaction $\ce{Li^+ + e^- <=> Li + \kappa S}$ we have $\lambda := \mu_{Li}|^\A + \kappa \mu_S|^\E - \mu_{Li^+}|^\E - \mu_{e^-}|^\A$ whereby $\lambda > 0$ entails $ r = L \cdot g(\lambda) > 0$. Since $J_\A = - \hat D_\A \nabla y_\A$ we have at the boundary 
\begin{align}
    J_\A \cdot \vecn =  (+1) r 
\end{align}
where (+1) is the stoichiometric coefficient of the product \ce{Li}. %
%
%

\subsection{Initial values and potential} 
\label{sub:inital_values}
The initial values, also used for the Newton solver, write with
\begin{align}
        y_\A^0 &:= \begin{cases}
                        y_\A^{\Cat,0}           & \qquad                (\xi,r) \in \Omega_\A^\Cat \\
                        y_\A^{\An,0}            &\qquad                 (\xi,r) \in \Omega_\A^\An \\            
        \end{cases} ~,&&& 
        y_\E^0 &:= \frac{n_\EC}{n_\ES - 2 \cdot \kappa_\E \cdot n_\EC} \\
        \tilde \phi_\S^0 &:= \begin{cases}
                f_\A^\An(y_\A^{\An,0}) - f_\A^\Cat(y_\A^{\Cat,0}))  & \qquad \xi \in \Omega^\Cat \\
                0                                                         &             \qquad   \xi \in \Omega^\An \\          
        \end{cases} ~, &&&
        \tilde \phi_\E^0 &:=f_\A^\An(y_\A^{\An,0}) - f_\E(y_\E^0) 
\end{align}
simply as 
\begin{align}
        y_\A(\xi,r,\tau=0) &= y_\A^0 & 
        y_\E(\xi,\tau=0) &=y_\E^0       \\
        \tilde\phi_\E(\xi,\tau=0) &=\tilde \phi_\E^0 &
        \tilde\phi_\S(\xi,\tau=0) &=\tilde \phi_\S^0 ~.
\end{align}
\subsubsection{Parameters} 
\label{ssub:parameters}
All parameters and their values for the subsequent numerical calculations are summarized in appendix \ref{app:par}. 

\section{Discretization and Model Order Reduction of the Battery Model}
\label{sec:headings}
The modeling approach discussed in the previous section neglects the details of the electrodes microstructure and describes it as a homogeneous medium in which electrolyte and the solid electrode materials exist at every point. This homogenized model is called macroscopic model in the following. In this macroscopic model, the intercalation of Li-ions in the electrode particles is incorporated through a coupled diffusion equation in radial direction of the particles in each macroscopic quadrature point.
In this way we get a pseudo-2D model (i.e. 1D+1D model) for the full battery cell. The model is given by a system of nonlinear PDEs for the homogenized electric potentials $(\tilde{\varphi}_E , \tilde{\varphi}_S )$ and the mole fraction of $\text{Li}^+$-ions $(y_E, y_A)$ in the electrolyte and in the positive and negative electrode materials, respectively.
The equation for the mole fraction in the electrode materials is calculated in the additional pseudo dimension, namely in the particle radius $\nu$. The overall PDE system can be written in the following abstract  form:
\begin{align*}
- \alpha(\cdot, y_A) \, \frac{\partial y_A}{\partial \tau}& - \nabla_{\nu} \cdot  [- \beta(\cdot,y_A) \nabla_{\nu} y_A] & =& \,\,  0, \\
&- \nabla_{\xi} \cdot [- \gamma(\cdot) \nabla_{\xi}  \tilde{\varphi}_S]  &- R_1(\cdot ,\tilde{\varphi}_E, \tilde{\varphi}_S, y_E,  y_A \vert_ {\nu=1}) =& \,\, 0 , \\
-\delta(\cdot, y_E) \, \frac{\partial y_E}{\partial \tau}  &- \nabla_{\xi} \cdot [ - \kappa(\cdot,y_E) \, \nabla_{\xi} y_E ] & + R_2(\cdot ,\tilde{\varphi}_E, \tilde{\varphi}_S, y_E,  y_A \vert_ {\nu=1}) =& \,\, 0, \\
&- \nabla_{\xi} \cdot [- \omega (\cdot, y_E) \, \nabla_{\xi} y_E - \rho(\cdot, y_E) \, \nabla_{\xi} \tilde{\varphi}_E ]\!\!\!\!\! &\!\!\!\!\!+ R_3(\cdot ,\tilde{\varphi}_E, \tilde{\varphi}_S, y_E,  y_A \vert_ {\nu=1}) =& \,\, 0.
\end{align*}
The linear and nonlinear coefficient functions $\alpha, \beta, \gamma, \delta, \kappa, \omega$ and $ \rho$ correspond to the representation from the equations (2.129)-(2.132) and depend on the domain for which the system is defined (anode, cathode and separator). $R_i , i = \{1, 2, 3, 4\}$ represents the reaction rate functions (2.127) in addition to the previous constants. The system is completed by the boundary conditions as well as interface conditions (2.133)-(2.134). Note that there are corresponding Neumann boundary conditions $\beta(\cdot,y_A) \nabla_{\nu} y_A \cdot n =  - R_4(\cdot,\tilde{\varphi}_E, \tilde{\varphi}_S, y_E,  y_A)$ at the boundary of the electrode particles, i.e. at $\nu = 1$, which couples the microscopic equation with the macroscopic equations.

\subsection{Discretization}
\label{sec:discr}

For the battery model in the abstract form above, let $\Omega_{1D}$ denote a computational one-dimensional domain in the macroscopic direction and $\Omega_{\nu} = \Omega^x_{\nu}$ the transverse/radial directions in the electrode particles associated with each $x \in \Omega_{1D}$. Let $\mathscr{P} \subset \mathbb{R}^P, P \geq 1$ denote the parameter space. We define the solution space $V = V_1 \oplus V_2$ with $ V_1 = H^1(\Omega_{\nu}), V_2 = (H^1(\Omega_{1D}))^3$ and $V^{\prime}$, the dual space of $V$. 
For a corresponding variational weak formulation, we obtain, after a semi-discretization in time $t$ by the implicit Euler method, that the battery model can be formulated as the following nonlinear system:
\begin{align}
\text{Find}\ u^{t+1}=[u^{t+1}_1, u^{t+1}_2, u^{t+1}_3, u^{t+1}_4]^\top \in V: \quad G_{\mu}(u^{t+1},\psi)=f_{\mu}(u^{t},\psi) \ \ \forall \ \psi \in V,
\label{G_fom}
\end{align}
where the operator $G_{\mu}(\cdot,\cdot): V \times V \rightarrow \mathbb{R}$ represents the non-linear time-discrete PDE system. The index $\mu \in \mathscr{P}$ indicates the dependence of the problem on certain parameters, such as the charge rate, the diffusion coefficient or the reaction rate. $f_{\mu}(u^{t}, \cdot) \in V^{\prime}$ contains the solid potential Neumann boundary conditions.

For the discretization in space, the finite element method is used \cite{FE}, 
i.e. we project (\ref{G_fom}) to a finite dimensional, continuous and piecewise polynomial space $V_h \subset V$. 
We hence obtain for each time step a fully-discrete non-linear system of the form:
\begin{align}
\text{Find}\ u_h^{t+1}=[u^{t+1}_{1_h}, u^{t+1}_{2_h}, u^{t+1}_{3_h}, u^{t+1}_{4_h}]^\top \in V_h: \quad G_{\mu}(u_h^{t+1},\psi_i)=f_{\mu}(u_h^{t},\psi_i) \quad \forall \, i=1,\dots n,
\label{GS_fom}
\end{align}
where $\psi_i, i=1, \dots, n$ denotes the standard Lagrange basis of the finite element space $V_h = \bigoplus_{i=1}^4 V_{i_h}$. Henceforth, the operator $G_{\mu}$ can be called the finite element operator which operates on $V_h \times V_h$.  The developed discretization does not depend on the specific choice of the parameters to be varied.


\subsection{Reduced Basis Method}
\label{sec:rb}
As e.g. detailed in  \cite{RB_Manzoni}, the reduced basis method is based on the idea of performing a Galerkin projection of the discrete equations onto low-dimensional subspaces $\tilde{V} \subset V_h$ in order to accelerate the repeated solution of (\ref{GS_fom}) for varying parameters $\mu$. Under this projection, the reduced problem is given by:
\begin{align}
\text{Find}\ \tilde{u}^{t+1}=[\tilde{u}^{t+1}_1, \tilde{u}^{t+1}_2, \tilde{u}^{t+1}_3, \tilde{u}^{t+1}_4]^\top \in \tilde{V}: \quad G_{\mu}(\tilde{u}^{t+1},\tilde{\psi}_i)=f_{\mu}(\tilde{u}^{t},\tilde{\psi}_i)\ \forall \, i=1,\dots m,
\label{GS_rom}
\end{align}
where $m \ll n$ and $\tilde{\psi}_i, i=1,\dots,m$, represents the basis of the reduced space $\tilde{V}$.
The basic idea of the reduced basis method is to perform a so-called offline/online decomposition.
In the preceding offline phase, the reduced space is constructed, and after projecting onto this reduced space, the resulting low-dimensional problem can be solved for any suitable parameter value in a following online phase. Various methods for constructing the reduced space for time-dependend problems have been considered in the literature, such as the POD-Greedy \cite{Pod_Greedy}. The POD-Greedy method produces approximation spaces with quasi-optimal $l^{\infty}$ -in-$\mu$, $l^2$ -in-time reduction error \cite{Pod_Greedy_RE}.

As the POD-Greedy method relies on the usage of rigorous and efficient to evaluate a posteriori error estimators, which are not available for the non-linear battery model at hand, in our numerical experiments, the reduced space is generated by the PODs of pre-selected set of solutions trajectories of the problem (\ref{GS_fom}), called snapshots.
Let $\mathcal{P} = \{ \mu^1, \dots, \mu^{n_s} \}$ be a set of $n_s$ parameter samples and $\{ u_h(\mu^1), \dots, u_h(\mu^{n_s}) \}$ the corresponding snapshot set. At each time step the snapshot of the set is calculated using the following prescription of Newton's method:
\begin{align*}
D_uF_{\mu}(u_h^{t+1,k},\psi_i) \ \delta u_h &= - F_{\mu}(u_h^{t+1,k},\psi_i) \qquad \forall \, i=1,\dots,n,\\
u_h^{t+1,k+1}(\mu) &= u_h^{t+1,k}(\mu) + \delta u_h,
\end{align*}
where $D_u F_{\mu}(z,\psi_i)$ is the Fr\'{e}chet derivative of $F_{\mu}           (u,\cdot) = G_{\mu}(u,\cdot) - f_{\mu}(u_h^t,\cdot)$ with respect to $u$ at $z \in V_h$.

We separate the snapshots into the respective components, $u_{1_h}, \ldots, u_{4_h}$ and generate the reduced basis separately, cf. \cite{Stephan_Felix}. We define the corresponding snapshot matrices $S_i \in \mathbb{R}^{N_h^{\iota} \cdot \vert t \vert \times n_s}$ with $i=1,2,3,4$ and $\iota \in \{ \Omega_{1D}, \Omega_{\nu}\}$ as:
\begin{align*}
S_i &= \left[ u_i^1, \dots, u_i^{n_s}\right],
\end{align*}
where the vectors $u_i^j \in \mathbb{R}^{N_h^{\iota} \cdot \vert t \vert},\, 1 \leq j \leq n_s,\,$ denote the degrees of freedom  of the functions $u_{i_h}(\mu^j)\vert_{t} \in V_{i_h}$. The singular value decomposition of $S_i$ is given through
$
S_i = U_i \Sigma_i Z_i^T,
$
where $U_i \in \mathbb{R}^{N_h^{\iota} \cdot \vert t \vert \times N_h^{\iota} \cdot \vert t \vert}$ and $Z_i \in \mathbb{R}^{ns \times ns}$ represent orthogonal matrices, and $\Sigma_i = \text{diag}(\sigma_i^1, \dots, \sigma_i^{z_i}) \in \mathbb{R}^{N_h^{\iota} \cdot \vert t \vert \times ns}$ with $\sigma_i^1 \geq \sigma_i^2 \geq \dots \geq \sigma_i^{z_i}, \, z_i \leq \min(N_h^{\iota} \cdot \vert t \vert, ns)$ contain the singular values. The left singular vectors 
\begin{align*}
U_i = \left[ \zeta_i^1 |\dots | \zeta_i^{N_h^{\iota}}  \right]
\end{align*}
span the reduced space $\tilde{V_{i}}$ using only the singular vectors whose singular values are above a fixed threshold value. Due to the fundamental properties of POD, the projection error consist of the $l^2$-sum of the corresponding truncated singular values, and $\tilde{V}_i, i=1,2,3,4,$ are $l^2$ -in-space, $l^2$ -in-time best-approximation spaces for the considered training set of snapshots. The overall reduced space is defined as $\tilde{V} := \bigoplus\limits_{i=1}^4 \tilde{V_{i}}$.

\subsection{Empirical Interpolation and Hierarchical Approximate POD}
\label{sec:eim}

The model reduction approach described in Subsection \ref{sec:rb} 
still suffers from a huge offline cost, due to the global POD construction of the reduced basis and from an inefficient online computational cost, as the reduced model (\ref{GS_rom}) cannot be decomposed efficiently into high- and low-dimensional computations due to the presence of non-linearities in the system. In order to obtain a more efficient ROM with reduced computational cost in the offline phase and a full so called offline-online decomposition we 
replace the global POD reduced basis construction by an incremental hierarchical approximate POD (HAPOD) \cite{HAPOD} and construct an affine approximation of the non-linear differential operator by an empirical operator interpolation \cite{deim, eim_operator}.

In detail, the offline-online decomposition
consists of precomputing parameter- and solution-independent terms, a so called collateral basis that allows to interpolate evaluations of the operator $G_{\mu}$. The construction of the collateral basis and associated interpolation functionals is done once in the offline phase to allow a fast evaluation of the interpolated operator $I_M(G_{\mu})$ during the online phase. in detail, the empirical operator interpolation generates a separable approximation by interpolating at $M$ selected degrees of freedom of $V_h$. The approximation of the operator is of the following form:
\begin{align}
G_{\mu}(\tilde{u},\cdot) \approx I_M(G_{\mu}(\tilde{u},\cdot)) = \sum_{q=1}^M \theta_{\mu}^q(\tilde{u}) G^q,
\end{align}
where $\{ G^q\}_{q=1}^M$ is the collateral basis, i.e. a basis for a subspace of $$\mathscr{M}_G = \{ G_{\mu}(u_h,\cdot) | \mu \in \mathscr{P} \},$$
and $\theta_{\mu}^q$ are interpolation coefficients recalculated for each $\mu$ and $\tilde{u}$ during the online phase. The collateral basis is obtained by applying the POD method to the set of snapshots obtained as images under the operator, i.e.  $G_{\mu}(u_h,\cdot)$. The set of  snapshots  thereby  includes the Newton stages in addition to the corresponding solution trajectory of $G_{\mu}(u_h,\cdot)$. Moreover, in analogy to the construction of the reduced basis described above, also the collateral basis is constructed  separately for each solution component. Hence, we define for $i=1,2,3,4$:
\begin{align*}
\left[ G_i^1, \dots, G_i^{M_i} \right] = \text{POD}\left( [G_{\mu_1}^1(u_i), \dots, G_{\mu_1}^{s_1}(u_i), G_{\mu_2}^1(u_i), \dots, G_{\mu_{n_s}}^{s_{n_s}}(u_i)], \epsilon_{POD} \right),
\end{align*}
where the vectors $G_{\mu_j}^k(u_i) \in \mathbb{R}^{N_h^\iota}, 1 \leq j \leq n_s,\, 1 \leq k \leq s_j,$ represent the degrees of freedom  of the functions $G_{\mu_j}(u_{i_h}^k(\mu^j))$ and $\epsilon_{POD}$ the error tolerance for the POD. We define $\left[ G^1, \dots, G^M \right] = \left[ G_1^1, \dots, G_4^{M_4} \right]$ with $M = \sum_{i=1}^4 M_i$.

To calculate the interpolation coefficients $\theta_{\mu}^q(\tilde{u}), q=1, \dots,M$ for given $\mu$ and $\tilde{u}$, the interpolation constraints are imposed at $M$ interpolation points. The interpolation points are selected iteratively from the indices of basis $\{ G^q\}_{q=1}^M$ using a greedy procedure.
This procedure determines each new interpolation point by the minimization of the interpolation error over the snapshots set measured in the maximum norm. For more details we refer to \cite{deim, eim}. 

By replacing the operator $G_{\mu}$ in (\ref{GS_rom}) by the fast to evaluate interpolated operator $I_M (G_{\mu})$, we obtain the completely offline-online decomposable reduced problem
\begin{align}
\text{Find}\ \tilde{u}^{t+1} \in \tilde{V}: \quad I_M (G_{\mu}(\tilde{u}^{t+1},\tilde{\psi}_i))=f_{\mu}(\tilde{u}^t,\tilde{\psi}_i)\ \forall \, i=1,\dots m,
\label{GS_rom+EI}
\end{align}
which can be solved efficiently for varying parameters $\mu$. 

For large-scale time dependent applications such as our battery model, computing the POD algorithm can be expensive. Especially if we include the evaluations of $G_{\mu}$ at all Newton levels of the selected solution trajectories in the operator snapshot set. The hierarchical approximate POD (HAPOD) algorithm is an efficient approach, which approximates the standard POD algorithm based on tree hierarchies, where the task of computing a POD for a given large snapshot set $S$ is replaced by multiple small PODs \cite{HAPOD}. More specifically, we use the special case of incremental HAPOD. In this case, the tree structure is such that each node of this tree represents either a leaf or has exactly one leaf and one non-leaf as children.  
In detail, first the vectors of the given snapshot set are assigned to the leaves of the tree $\beta_1, \cdots , \beta_B$. Starting with two leaves $\beta_1, \beta_2$, a POD of each local snapshot data is computed. The resulting modes scaled by their singular values are the input to the parent node $\alpha_1$, which is again used to calculate a POD. This newly generated input and the local snapshot data assigned to leaf $\beta_3$ are the input of the parent node $\alpha_2$. The final HAPOD modes are reached, when the last leaf $\beta_B$ has entered. For the calculation of the collateral basis, we e.g. define for $i = 1,2,3,4$:
\begin{align*}
\left[ G_i^1, \dots, G_i^{M_i} \right] = \text{HAPOD}\left( [G_{\mu_1}^1(u_i), \dots, G_{\mu_1}^{s_1}(u_i), G_{\mu_2}^1(u_i), \dots, G_{\mu_{n_s}}^{s_{n_s}}(u_i)], \epsilon_{POD}, \omega \right),
\end{align*}
where $\epsilon_{POD}$ 
is the desired approximation error tolerance for the resulting HAPOD space. Depending on omega, one might get more modes than needed for a POD with the same tolerance. The local tolerances in the HAPOD algorithm are computed from $\epsilon_{POD}$ and $\omega \in (0,1)$. More details can be found in \cite{HAPOD}.

\section{Numerical Results}
\label{sec:results}
In order to create a test environment for our modeling framework, we developed an experimental implementation of the aging effects of the battery model from Section \ref{sect:derivation_model}. We will investigate the efficiency of the reduced order simulations by evaluating electrochemical characteristics over the cyclization $n = 1, \dots, 1000$ for different aging models. The electrochemical characteristics are the voltage-capacity spectrum $E(\bar y_A , n; C_h)$ and the status of charge $\bar y_A^{Cat}$ at a specific voltage value $E_{min}$ (see eq. (2.102)-(2.109)).
We assume that the aging effects are modeled by given functions in dependence of the cycle number $n$ for the reaction rate $\hat L(n)$ and diffusion coefficient $\hat D_A(n)$. These functions are used to investigate the qualitative behavior of the aging effects. The parameter dependence of the reaction rate $\hat L$ examines the degradation of the solid electrolyte interphase. This illustrates the increase in reaction resistance due to cyclization.
Furthermore, the degradation of the porous electrodes is investigated and represented by a decrease in the diffusion coefficient $\hat D_A$. This effect is caused by micro-cracks within a particle.

To efficiently analyze this forward modeling of aging effects, we consider three scenarios. In the first scenario, we consider the unaged battery and calculate the voltage against the state of charge for varying charge rates $C_h$. In the next case, we set $C_h = 1$ and examine the aging effects of $\hat D_A$ and $\hat L$ by alternately choosing one of the parameters fixed. In the last case we vary all parameters $\hat D_A, \hat L$ and $C_h$.

\subsection{Implementation Aspects}
\label{sec:impl}
\begin{figure}[t!] 
    \centering
    \subfloat{\includegraphics[height=5cm, width=0.45\textwidth]{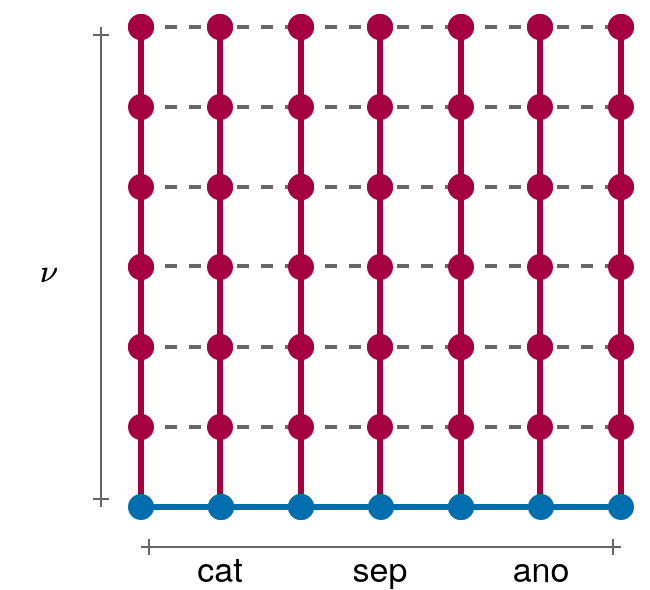}}
    \caption{Sketch of the pseudo 2D grid. The blue line shows the computational domain of the macroscopic equations, while the red lines illustrate the computational domain for the microscopic equation.}
\label{2Dgrid}
\end{figure}
Let us first introduce the settings used in the numerical experiments. We implemented a (pseudo) 2D grid, where the bottom of the grid corresponds to the 1D grid on which the macroscopic equations are computed. The length of the bottom $L_{cat} + L_{sep} + L_{ano}$ is divided into $N_h^{\Omega_{1D}}=300$ grid points. Here, each of the cathode, separator and anode is discretized into 100 grid points. The pseudo-dimension for the microscopic equation, which goes vertically at each bottom grid point, consists of $N_h^{\Omega_{\nu}}$ grid points. In the simulations $N_h^{\Omega_{\nu}} = 100$ was chosen. All in all, this results in 30900 degrees of freedom (dofs) for the overall system. Thereby, we assume that the high-dimensional discretization adopts a resolution of $\vert\vert u_n - u_m \vert\vert_2 \leq 10^{-5}$ with 30900 dofs for $u_n$ and 31512 dofs ($N_h^{\Omega_{1D}}=303, N_h^{\Omega_{\nu}} = 101$) for $u_m$. The Neumann boundary condition of the microscopic equation couples the equation to the macroscopic equations and is defined for $\nu = 1$. To ensure that the bottom grid corresponds to $\nu=1$ a transformation of the microscopic equation with $\tilde{\nu}=(1-\nu)$ is performed. In addition, as an essential step for the stability of the model, a variable transformation of the microscopic variable $y_A$ of the following form is applied: 
\begin{align*}
y_A(g) = \frac{e^g}{1+ e^g}, \quad g(y_A) = \text{ln} \left( \frac{y_A}{1-y_A} \right).
\end{align*}
The time discretization is performed by an implicit Euler method on a $T = 1$ time interval with a time step size of $\Delta t = 10^{-2}$.
The nonlinear battery system is solved by a Newton method to a relative error accuracy of $10^{-5}$  and a termination condition of $\min_{x \in \Omega^{Cat}} \tilde{\varphi}_S(x) \le E_{min}$ with the voltage value $E_{min}=-0.2 $.

To generate the reduced space $\tilde{V}$, we compute a snapshot set $S_{train}$ on training sets of equidistant parameters. For the experiments 1 and 3 we choose $\vert \mathcal{P}_{train} \vert = 15$ and for experiment 2 we choose $\vert\mathcal{P}_{train} \vert= 10$. 
As a measure for the model reduction error we determine the relative $l^2$ -in-space, $l^2$ - in-time error averaged over a set of random test parameters $\mathcal{P}_{test}$ given by: 
\begin{align}
 \frac{1}{\vert \mathcal{P}_{test} \vert}\sum_{\mu \in \mathcal{P}_{test}} \frac{\vert\vert u_h(\mu) - \tilde{u}(\mu) \vert\vert_2}{\vert \vert \tilde{u}_(\mu) \vert \vert_2}.
 \label{l2l2error}
\end{align} 
The programming language is Python. All simulations of the high-dimensional model are computed with the finite element sofware NGSolve \cite{ngsolve}. NGSolve provides the ability to construct the complex grid structure and define the battery model for each subdomain.
For the implementation of the reduced basis method, the NGSolve code has integrated into the model order reduction library pyMOR \cite{pymor}.
We include the evaluations of $G_{\mu}$ on all Newton stages of the selected solution trajectories in the operator snapshot set to compute the collateral basis. This leads to a stabilization of the reduced model.
In order to speed up the computation of the collateral basis for the empirical interpolation data via POD, the HAPOD algorithm is used instead of the standard POD algorithm. For the generation of the reduced basis, we use the HAPOD algorithm as well. We choose $\epsilon = 4e-8$ and $\omega = 0.9$ in both cases.
For illustration purposes, the reduced space and empirical interpolation are calculated for a training set $\vert\mathcal{P}_{train} \vert= 5$ for the variation of $C_h$. In this case, the calculation without using HAPOD takes $107.32$ min, while only $7.31$ min CPU-time were needed using HAPOD. This corresponds to a speedup of $14.68$ in the offline phase. All tests were performed on the same computer and software basis.

\subsection{Experiment 1}
\label{sec:exp1}
\begin{table}
\caption{Relative model reduction error (\ref{l2l2error}) and reduced simulation times for a battery simulation trajectory with $\hat L=0.5, \hat D_A=0.5$ and $\mathcal{P}_{test} = 10$ . The number of the reduced basis consists of the four variables, e.g. $11 = \# u_1 + \#u_2 + \#u_3 +\#u_4 = 2+2+4+3$. When the reduced basis is increased, each variable is added one basis. The number of interpolation point is 102. The average time for the solution trajectory of the high-dimensional model is 363.48 s.} \label{table1}
\begin{minipage}{\textwidth}
\tabcolsep=8pt
\begin{tabular}{lllll}
\hline\hline
{reduced basis size}  & {11} & {15} & {19} & {23} \\
\hline
relative error  & $2.39 \cdot 10^{-4}$ & $3.23 \cdot 10^{-5}$ & $1.62 \cdot 10^{-5}$ & $1.19 \cdot 10^{-5}$\\
time (s) & $25.85$ & $26.51$ & $26.94$  & $27.38$\\
speed up  & $14.11$ & $13.66$ & $13.55$ & $13.06$\\
\hline\hline
\end{tabular}
\end{minipage}
\end{table}
\begin{figure}[h] 
    \centering
    \subfloat[][]{\includegraphics[height=4.5cm, width=0.44\textwidth]{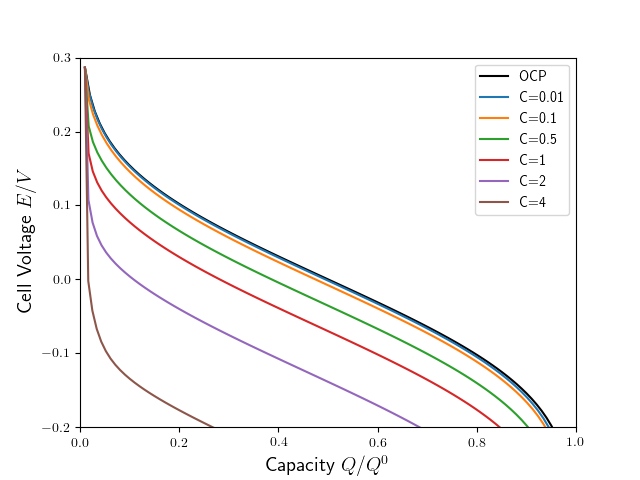}} \\ \vspace{-0.3cm}
    \subfloat[][]{\includegraphics[height=5.5cm, width=0.55\textwidth]{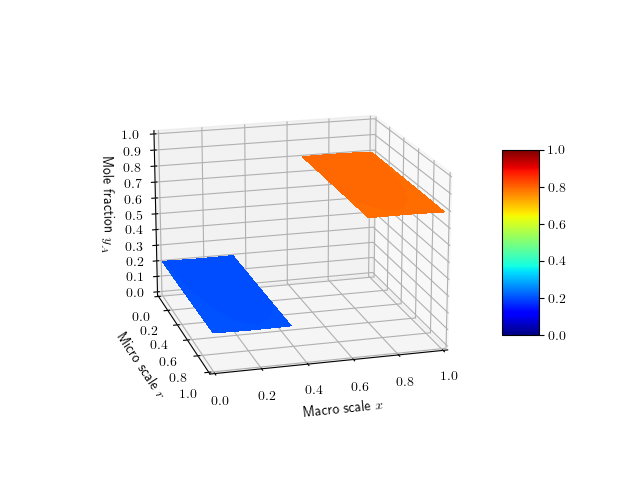}}
    \subfloat[][]{\includegraphics[height=4.5cm, width=0.44\textwidth]{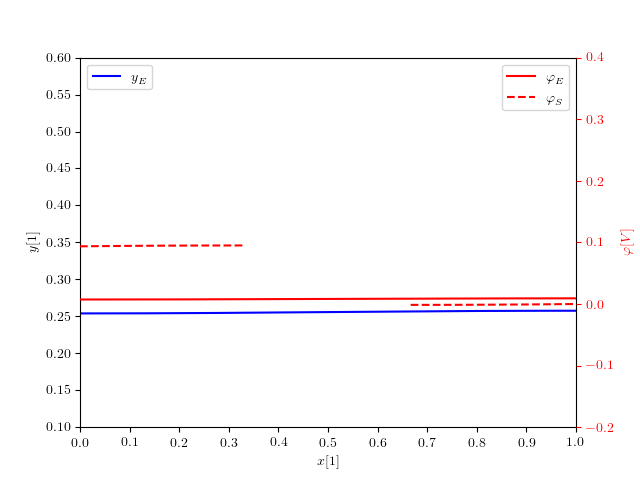}} \\ \vspace{-0.3cm}
    \subfloat[][]{\includegraphics[height=5.5cm, width=0.55\textwidth]{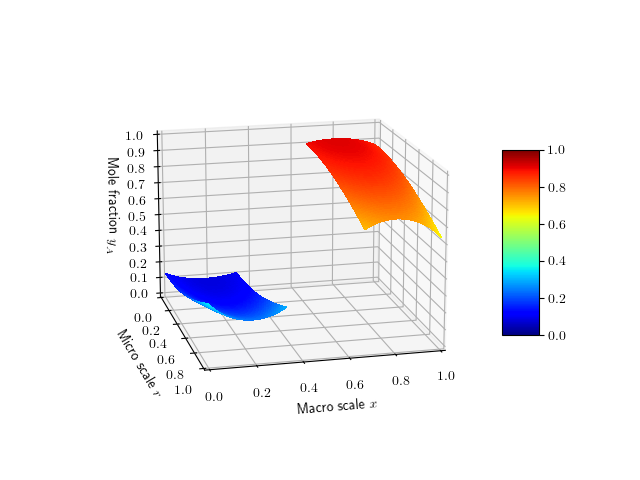}}
    \subfloat[][]{\includegraphics[height=4.5cm, width=0.44\textwidth]{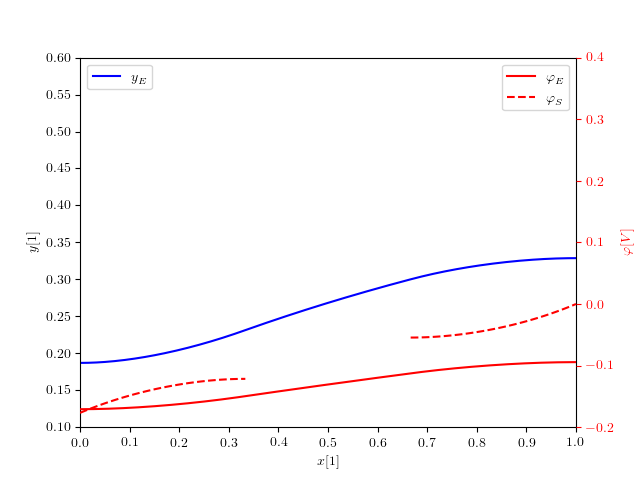}}
    \caption{a) Voltage-capacity spectrum compared to the open circuit potential. Solution plots of the four components, b)-c) with $C_h = 0.1$ and d)-e) with $C_h = 4$ for $t=0.2$.}
\label{exp1}
\end{figure}

In this subsection, we consider the variation of the charge rate $C_h$ with $C_h \in [0.01,4]$ for an unaged battery. In this case, we choose $\hat L=0.5$ and $\hat D_A=0.5$. For the reduced order model the number of basis functions for the four variables are set to $ \# u_1  = 3,  \#u_2 = 3, \#u_3 = 5$ and $\#u_4 = 4$. For the empirical operator interpolation, the number of interpolation points are $ \# G(u_1)  = 19,  \#G(u_2) = 15, \#G(u_3) = 60$ and $\#G(u_4) = 8$. These numbers are obtained by calculating the relative model reduction error (\ref{l2l2error}) with successive increase of the basis size up to an accuracy of order $10^{-5}$ (see Table \ref{table1}). To ensure the stability of the reduced model for empirical operator interpolation, the number of interpolation points must be chosen large enough. Especially the number of interpolation points for $G(u_3)$ are crucial here. \\
The voltage-capacity spectrum is shown in Fig. \ref{exp1}a), where we achieve a model reduction error of less than $10^{-4}$ for a simulation time of 2.58 minutes. Therefore, we obtain a speed up of 15.41.
Two solution plots at a fixed time $t=0.2$ for the charge rate $C_h \in \{ 0.1, 4 \}$ are illustrated in Fig. \ref{exp1}b)-e). At a low charge rate, we almost reach the open circuit potential (OCP), which can be observed by the fact that nearly constant functions are obtained. For larger charge rates, we observe higher gradients in macroscopic and microscopic directions due to transport limitations.

\subsection{Experiment 2}
\label{sec:exp2}
\begin{table}
\caption{Relative model reduction error (\ref{l2l2error}) and reduced simulation times for a battery simulation trajectory with $C_h=1$ and $\mathcal{P}_{test} = 10$. The number of the reduced basis consists of the four variables, e.g. $10 = \# u_1 + \#u_2 + \#u_3 +\#u_4 = 2+2+4+2$. In each column, a basis is added to each variable. The number ob interpolation points amounts to 42. The average time for the solution trajectory of the high-dimensional model is 356.49 s.} \label{table_D+L}
\begin{minipage}{\textwidth}
\tabcolsep=8pt
\begin{tabular}{lllll}
\hline\hline
{reduced basis size} & {10} & {14} & {18} & {22}\\
\hline
relative error & $1.65 \cdot 10^{-5}$ & $6.43 \cdot 10^{-6}$ & $3.89 \cdot 10^{-6}$ & $ 1.30 \cdot 10^{-6}$\\
time (s) & $8.63$   & $8.42$ & $8.66$ & $8.83$\\
speed up & $41.27$ & $42.38$ & $41.07$ & $40.26$\\
\hline\hline
\end{tabular}
\end{minipage}
\end{table}
\begin{figure}[t!] 
    \centering
    \subfloat{\includegraphics[height=5cm, width=0.5\textwidth]{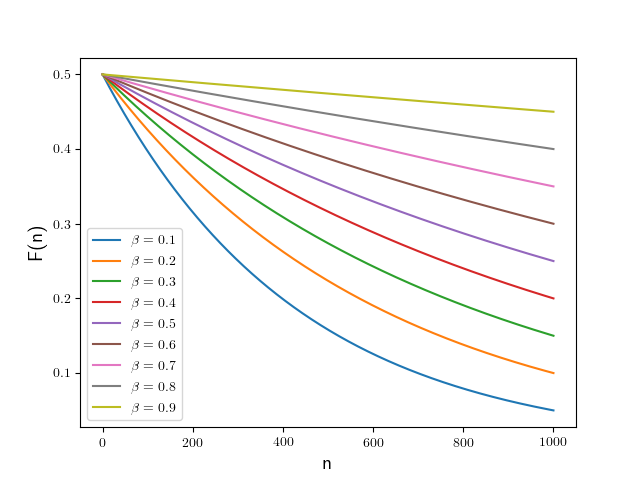}}
    \caption{Various evolutions of the parameter functions satisfying the partial differential equation \ref{ode_D+L} with $F_0 =0.5$ and $N = 1000$. }
\label{exp2_beta}
\end{figure}
\begin{figure}[h!] 
    \centering
    \subfloat[][]{\includegraphics[height=4.5cm, width=0.44\textwidth]{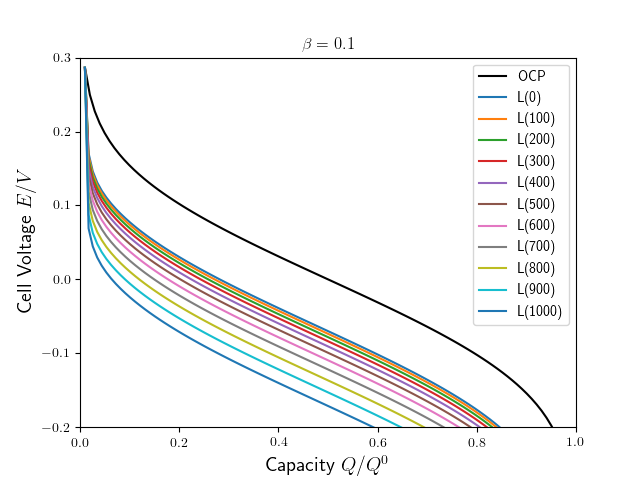}} 
     \subfloat[][]{\includegraphics[height=4.5cm, width=0.44\textwidth]{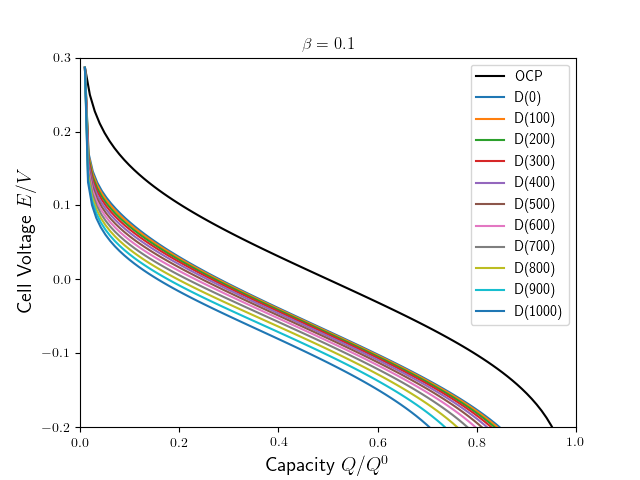}}  \\ 
    \subfloat[][]{\includegraphics[height=4.5cm, width=0.44\textwidth]{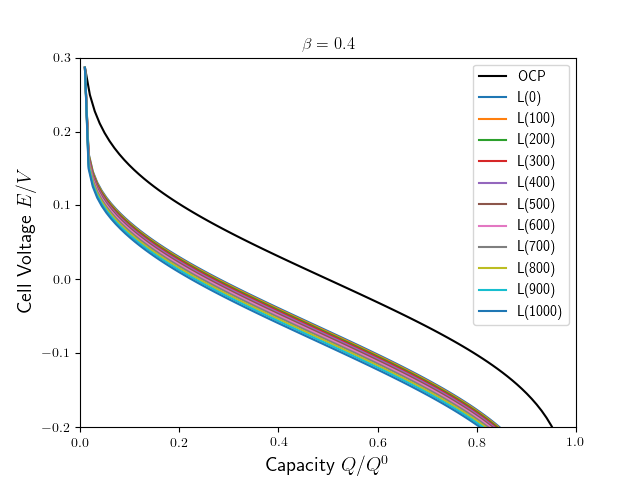}} 
     \subfloat[][]{\includegraphics[height=4.5cm, width=0.44\textwidth]{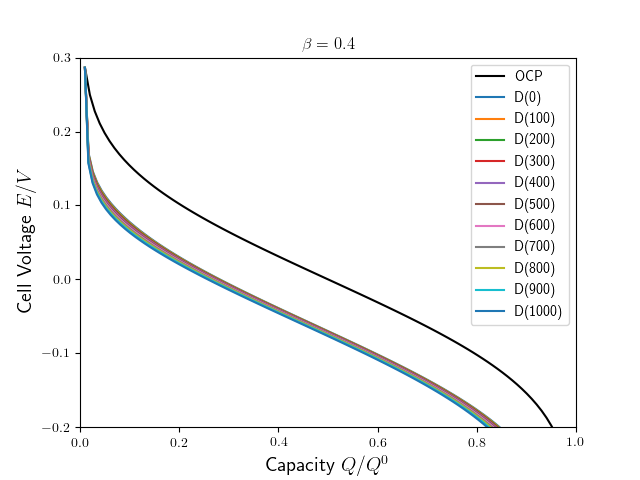}}  \\
    \subfloat[][]{\includegraphics[height=4.5cm, width=0.44\textwidth]{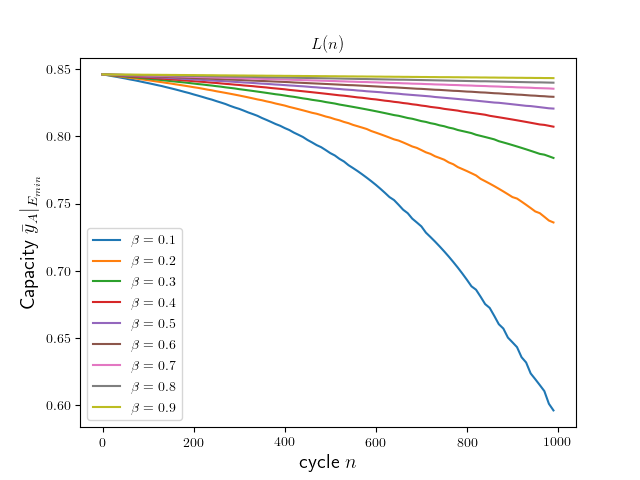}}
    \subfloat[][]{\includegraphics[height=4.5cm, width=0.44\textwidth]{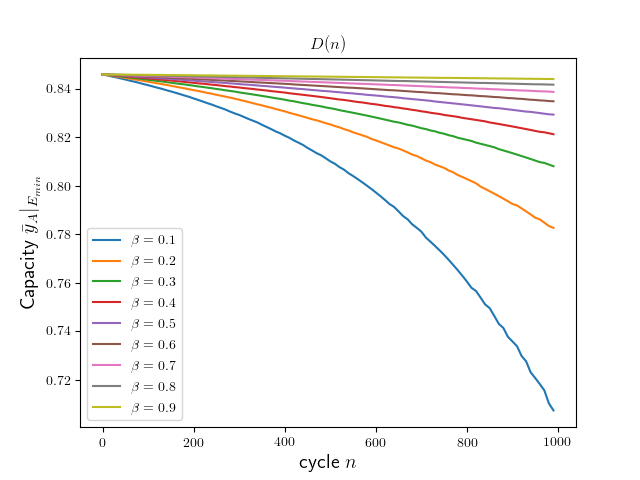}}\\
    \caption{a)-d) Evolution of the capacity-dependent voltage E of $\hat D_A(n)$ and $\hat L(n)$ compared to the open circuit potential for $\beta = 0.1$ or $\beta = 0.4$. e)-f) Effect of the different degradation models of $\hat D_A(n)$ and $\hat L(n)$ on the capacity at voltage $E_{min}$ over the number of cycles $n$ at $C_h=1$.}
\label{exp2}
\end{figure}

In the following the evaluation of the degradation simulations of the solid electrolyte interphase and the porous electrodes are presented. To investigate the qualitative behavior of these  aging effects, we consider the decrease in value of the parameters $\hat L$ and $\hat D_A$ equally in anode and cathode over $n$ with $n = 0 ,\dots, N$ cycles and set $C_h =1$. A cycle consists of a discharge process, assuming that the battery is charged in such a way that the chosen initial conditions apply at the beginning of each cycle. In addition, we assume that the evolution of the parameters $\hat L(n)$ and $\hat D_A(n)$ as a function of the number of cycles $n$ satisfies an ordinary differential equation
\begin{align}
\frac{\mathrm{d} F(n)}{\mathrm{d} n} &= a_F \, F(n), \label{ode_D+L}
\end{align}
with the unknown parameter $a_F(\beta)$ such that:
\begin{align}
F(0) &= F_0,  \\ 
F(N) &= \beta F_0, \quad \beta < 1. \notag
\end{align}
It follows that $F(n)= F_0 \, e^{\nicefrac{ \log (\beta) n}{N}}$ and $a_F(\beta) = \frac{\log(\beta)}{N}$. $F_0$ is the corresponding initial parameter value. In our case, $\hat D_{A,0}= \hat L_0 = 0.5$.  Under this assumption, we impose that the aging mechanism in one cycle depends on the value of the degradation of the previous cycle. (This assumption is based on the fact that under laboratory conditions, the cell is always discharged in the same way.) 

The characteristic spectrum of cell voltage $E$ for the degradation of reaction rate $\hat L$ for $\hat D_A=0.5$ is shown in Fig. \ref{exp2}a),c) for two different choices of $\beta$. Furthermore, the capacity $\bar y_A^{Cat}$ at the specific voltage value $E_{min}=-0.2$ over the number of cycles $N = 1000$ with a variation of $\beta$ is illustrated in Fig. \ref{exp2}e). The same scenario is shown for the degradation of the diffusion coefficient $\hat D_A$ for $\hat L=0.5$ in Fig. \ref{exp2}b),d),f). As expected, the graphs show when the parameters degrade faster and more severely, the cell voltage and capacity decrease more rapidly.  

For the implementation of the reduced order model with dependence on the parameters $\mu=[ \hat D_A$,$\hat L ]$ the number of basis functions for the four variables are set to $ \# u_1  = 3,  \#u_2 = 3, \#u_3 = 5$ and $\#u_4 = 3$. For the empirical operator interpolation, the number of interpolation points are $ \# G(u_1)  = 9,  \#G(u_2) = 9, \#G(u_3) = 15$ and $\#G(u_4) = 9$. As in Experiment 1, these numbers are obtained by calculating the relative model reduction error (\ref{l2l2error}), taking into account an accuracy of order $10^{-6}$ with successive increase of the basis size. To ensure the stability of the reduced model for empirical operator interpolation, the number of interpolation points must be increased for larger reduced basis space dimensions.

In Table \ref{table_D+L} we observe a rapid decay of the model reduction error, which stagnates already for relatively small reduced basis space dimensions. In this manner, we obtain relative reduction errors as small as $10^{-5}$ with simulation times of less than 10 seconds.
When calculating the voltage spectra, we achieve an average relative reduction error of about $10^{-3}$ and an average speed up of 24.37. Calculating the capacity at the voltage value $E_{min}$ over the number of cycles requires about 118.92 hours for the full model. By using the reduced model, we obtain approximations in about 2.53 hours with a relative error of $10^{-5}$. It is a speed up of 46.83.

\subsection{Experiment 3}
\label{sec:exp3}
\begin{table}
\caption{Relative model reduction error (\ref{l2l2error}) and reduced simulation times for a battery simulation trajectory with $\mathcal{P}_{test} = 10$ . The number of the reduced basis consists of the four variables, e.g. $13 = \# u_1 + \#u_2 + \#u_3 +\#u_4 = 3+3+4+3$. When the reduced basis is increased, each variable is added one basis. The number ob interpolation points amounts is 160. The average time for the solution trajectory of the high-dimensional model is 357.31 s.} \label{table_D+L+C}
\begin{minipage}{\textwidth}
\tabcolsep=8pt
\begin{tabular}{lllll}
\hline\hline
{reduced basis size} & {13} & {17} & {21} & {25}\\
\hline
relative error & $2.25 \cdot 10^{-5}$ & $1.40 \cdot 10^{-5}$ & $1.31 \cdot 10^{-5}$ & $1.20 \cdot 10^{-5}$\\
time (s) & $37.34$   & $37.68$ & $39.51$ & $39.49$\\
speed up & $9.82$ & $9.36$ & $9.05$ & $8.92$\\
\hline\hline
\end{tabular}
\end{minipage}
\end{table}
\begin{figure}[h] 
    \centering
    \subfloat{\includegraphics[height=5cm, width=0.5\textwidth]{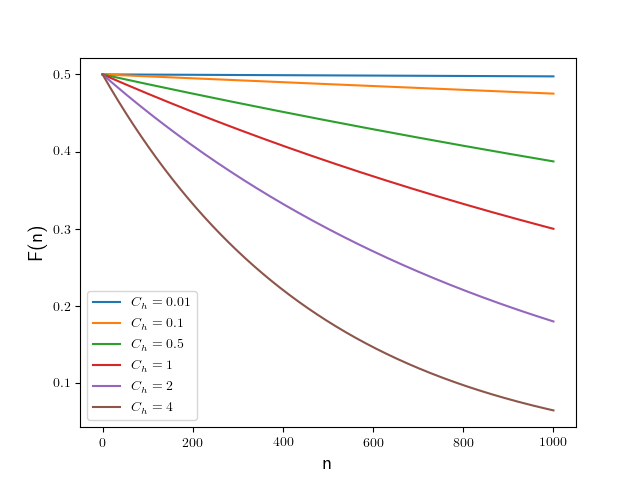}}
    \caption{Various degradation models that satisfy the partial differential equation \ref{ode_D+L+C} with $F_0 =0.5, \beta=0.6$ and $N = 1000$. }
\label{exp3_C}
\end{figure}

\begin{figure}[t!] 
    \centering
    \subfloat[][]{\includegraphics[height=4.5cm, width=0.44\textwidth]{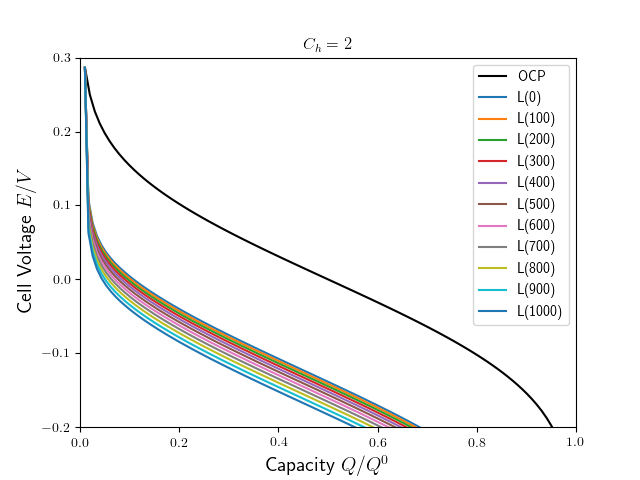}}
    \subfloat[][]{\includegraphics[height=4.5cm, width=0.44\textwidth]{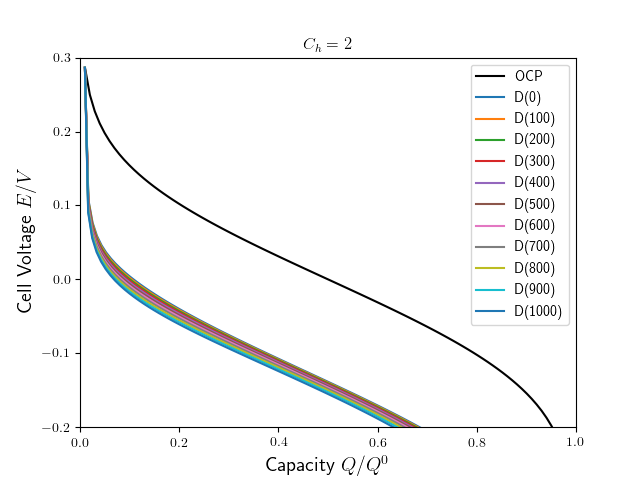}} \vspace{-0.3cm} \\ 
    \subfloat[][]{\includegraphics[height=4.5cm, width=0.44\textwidth]{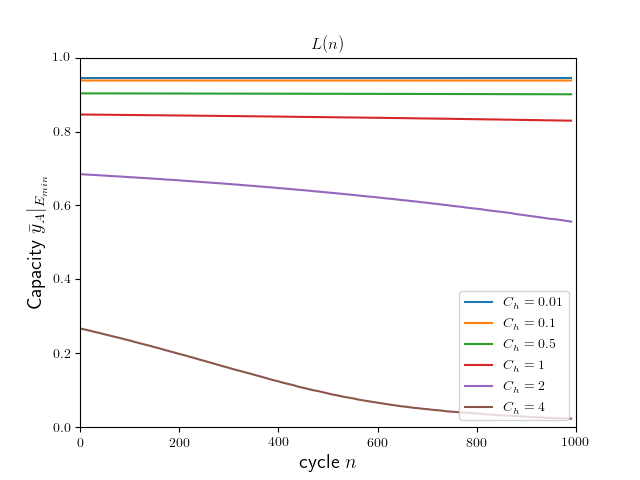}}
    \subfloat[][]{\includegraphics[height=4.5cm, width=0.44\textwidth]{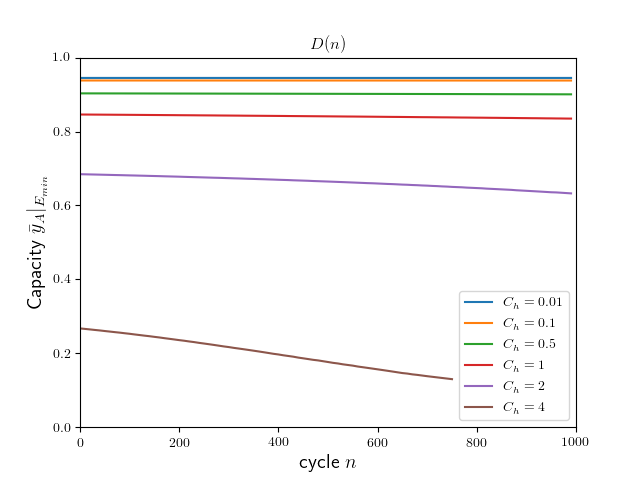}}\\
    \caption{a)-b) shows the spectrum of cell voltage E for the degradation of $\hat D_A(n)$ and $\hat L(n)$ compared to the open circuit potential for $C_h = 2$ . c)-d) illustrate the effect of the different degradation models of $\hat D_A(n)$ and $\hat L(n)$ on the capacity at voltage $E_{min}$ over the number of cycles $n$ with $\beta=0.6$.}
\label{exp3}
\end{figure}

In the previous section, degradation simulations of the solid electrolyte interphase and the porous electrodes were considered by decreasing the value of the parameters $\hat L$ and $\hat D_A$ over $N$ cycles. We assumed that the value of the parameter depends on the value of the parameter in the previous cycle. Thereby, the charge rate $C_h$ was set constant to $1$ for each simulation.
In this experiment, a variation of the charge rate is a matter of interest as well. This implies that the reduced order model depends on the following parameter vector $\mu=[C_h, \hat D_A, \hat L ]$. Furthermore, as in Experiment 2, we assume that the evolution of the parameters $\hat D_A$ and $ \hat L$ as a function of the number of cycles $n$ satisfies the ordinary differential equation that additionally depends on the charge rate $C_h$
\begin{align}
\frac{\mathrm{d} F(n)}{\mathrm{d} n} &= a_F \, F(n) \, C_h, \label{ode_D+L+C}
\end{align}
with the unknown parameter $a_F(\beta)$ such that:
\begin{align}
F(0) &= F_0,  \\ 
F(N) &= \beta F_0, \quad \text{for} \, \beta < 1 \, \text{at} \ C_h=1. \notag
\end{align}
It follows that $F(n)= F_0 \, e^{\nicefrac{ C_h \log (\beta) n}{N}}$ and $a_F(\beta) = \frac{\log(\beta)}{N}$. Here we set $\beta = 0.6$ and $F_0$ again represents the corresponding initial parameter value. Note that for $C_h=1 $ the situation is the same as in Experiment 2. 

In this experiment, we set the number of the basis functions for the four variables to $ \# u_1  = 4,  \#u_2 = 4, \#u_3 = 5, \#u_4 = 5$ and the number of interpolation points to $ \# G(u_1)  = 22,  \#G(u_2) = 20, \#G(u_3) = 110, \#G(u_4) = 20$, compare Table \ref{table_D+L+C}, taking into account an accuracy of order $10^{-5}$. Note, as in the first experiment, that the number of interpolation points, in particular the number of interpolation points for $G(u_3)$, must be chosen large enough to ensure the stability of the reduced model.

Comparing the run times from the calculation of the spectrum of the cell voltage (see Fig. \ref{exp3} a)-b)) for the full and reduced models, we obtain a speed up of about $7.97$ with a relative reduction error of about $10^{-5}$. In addition, the degradation of the capacity at the voltage $E_{min}$ over $N$ cycles is shown in Fig. \ref{exp3} c)-d). This resulted in a speed up of about $23.43$ with a relative reduction error of order $10^{-3}$.

\section{Conclusion}
\label{sec:conclusion}

We developed a mathematical model framework for an intercalation battery, consisting of multi-phase porous electrodes, on the basis of non-equilibrium thermodynamics. The framework is very flexible and applicable to a wide range of materials, either for the active intercalation phases or the (liquid) electrolyte, by a stringent formulation of the entire PDE problem in terms of general chemical potential functions. Special emphasis is put on thermodynamic consistency of the transport equations and their respective reaction boundary conditions by employing the very same chemical potential function entirely throughout the model. Periodic homogenization theory is applied to derive a general set of PDEs for the porous battery cell, where a special scaling of the micro-scale diffusion coefficient leads to a coupled micro-macro scale problem. Spherical symmetry of the intercalation particles is further employed, as well as a 1-D approximation of the macro-scale yields an effective 1D+1D non-linear PDE system. The (dis-)charge current, effectively characterized by the C-rate $C_h$, enters the PDE system as boundary condition for the electron flux. This allows, on the basis of numerical simulations, the computation of the time and space dependent thermodynamic state variables, \ie the electrolyte potential $\phi_\E(x,t)$, solid potential $\phi_\S(x,t)$, electrolyte concentration $y_\E(x,t)$, and active phase concentration $y_\A(x,r,t)$. Subsequently, this yields important characteristics of an intercalation battery, \ie the cell voltage $E$ as function of the status of charge $\bar y_\A$, parametrically dependent on the C-rate $C_h$. Further, battery degradation is considered in terms of cycle number $n$ dependent parameters, where exemplarily some degradation models in terms of simple evolution equations were stated. In order to simulate degradation effects, repeated numerical computations of the PDE system are required. For efficient numerical simulations, model reduction techniques were applied to the electrochemical battery model, \ie the reduced basis method combined with an empirical operator interpolation. We demonstrated the efficient applicability of these method with numerical studies on several aging scenarios. For degradation effects that impact the diffusion coefficient in the active phase or the intercalation reaction rate, we obtained capacity curves over the number of cycles with a speedup of about $46$, compared to to full numerical simulations of the same implementation. A speedup factor of about 23 was achieved by additionally investigating the effect of different choices of the charge rate. Numerical relative accuracy of order $10^{-3}$ (at least) was ensured within our simulations.

\appendix
\section{Parameters}
\label{app:par}

The following table summarizes all parameters and their values of the model in section \ref{sect:derivation_model}. \\

\begin{longtable}{|p{7cm}|p{4cm}|l|}
	\hline
	\textbf{Description}   & \textbf{Symbol} and \textbf{Value}   & \textbf{Units}  \\
	\hline
	\hline
	\textbf{Electrolyte} &&\\  
	\hline
	pure solvent molar concentration  	&  $n_\ESref = 11.9103 $						& $\uunits{\Umol \, \Ul^{-1}}$\\[0.1cm]
	reference electrolyte concentration &  $n_\ECref = 1 $								& $\uunits{\Umol \, \Ul^{-1}}$\\[0.1cm]
	solvent molar volume  				&  $v_\ES = \frac{1}{n_\ESref}$					& $\uunits{\Umol^{-1} \, \Ul}$\\[0.1cm]
	solvation number  					&  $\kappa_\EA = \kappa_\EC = 4$				& $\uunits{1}$\\[0.1cm]
	ion molar volume  					&  $v_\EA = v_\EC = \kappa_\EC \cdot v_\ES$		& $\uunits{1}$\\[0.1cm]
	molar conductivity 					&  $\tilde \Lambda_\E = 10$         	        & $\uunits{1}$ \\[0.1cm]
	transference number 				&  $t_\E = 0.5$									& $\uunits{\Uone}$ \\[0.1cm] 
	chemical diffusion coefficient    	&  $\tilde D_\E = 5$                            & $\uunits{1}$ \\[0.1cm]	
	\hline 
	\textbf{Cathode} (ideal electrode) &&\\
	\hline
	molar lattice concentration 		&  $n_{\A,\ell}^\Cat =   37.3114$				&   $\uunits{\Umol \, \Ul^{-1}}$ \\[0.1cm]
	initial value 		&  $y_\A^{\Cat,0} = 0.01 $				&   $\uunits{\Uone}$ \\[0.1cm]
	enthalpy parameter   				&  $\gamma_\A^\Cat = 1 $					&  	$\uunits{1}$ \\[0.1cm]
	electronic conductivity  			&  $\tilde \sigma_\S^\Cat = 10$		& 	$\uunits{\Uone}$\\[0.1cm]
	Li diffusion coefficient  			&  $\tilde D_{\A,\rref}^\Cat = 1$		& 	$\uunits{\Uone}$ \\[0.1cm]
	\hline
	\textbf{Anode}  (ideal electrode) &&\\
	\hline
	molar lattice concentration 		&  $n_{\A,\ell}^\An=  37.3114$			&   $\uunits{\Umol \, \Ul^{-1}}$ \\[0.1cm]
	initial value 		&  $y_\A^{\An,0} = 0.99 $				&   $\uunits{\Uone}$ 		\\[0.1cm]
	enthalpy parameter   				&  $\gamma_\A^\An = 1$	  				&  	$\uunits{1}$ \\[0.1cm]
	electronic conductivity  			& $\tilde \sigma_\S^\An = 10$		& 	$\uunits{\Uone}$ \\[0.1cm]
	Li diffusion coefficient  			&  $\tilde D_{\A,\rref}^\An = 1$		& 	$\uunits{\Uone}$ \\[0.1cm]
	\hline
	\textbf{Cathode Intercalation Reaction} &&\\  
	\hline 
	Half-Cell Reaction energy vs. metallic Li & $ E_{\A,\ce{Li+}}^\Cat = 3.95  $ & $\uunits{\UV}	$ \\[0.1cm]
	Exchange current density				  & $\tilde L^\Cat = 1$						& $\uunits{\Uone}$ \\[0.1cm]
	\hline
	\textbf{Anode Intercalation Reaction} &&\\  
	\hline 
	Half-Cell Reaction energy vs. metallic Li & $ E_{\A,\ce{Li+}}^\An=  0.2 $ & $\uunits{\UV}	$ \\[0.1cm]
	Exchange current density				  & $\tilde L^\An = 1$						& $\uunits{\Uone}$ \\[0.1cm]
	\hline
	\textbf{Geometry}   &&\\  
	\hline
	cathode thickness 	 &	$W^\Cat = 100 $ & $\uunits{\mu\Um}$\\
	separator thickness  &	$W^\Sep = 100$  & $\uunits{\mu\Um}$\\
	anode thickness    	 &	$W^\An  = 100$  & $\uunits{\mu\Um}$\\  
	\textbf{Micro-geometry}  &&\\  
    micro-unit cell width 			&  $\ell^\Cat = 10$ & $\uunits{\Unm}$ \\
	cathode particle radius 			&  $\tilde r_\A^\Cat = 0.4$	 & $\uunits{\Uone}$ \\
	cathode electrolyte phase fraction  &  $\psi_\E^\Cat  = 0.72713951$& $\uunits{\Uone}$\\[0.1cm]
	cathode electrolyte porosity tensor &  $\tenpi_\E^\Cat = 0.86842790$& $\uunits{\Uone}$   \\[0.1cm]
	cathode solid phase fraction   		&  $\psi_\S^\Cat = 0.27286022$ & $\uunits{\Uone}$ \\[0.1cm]
	cathode solid porosity tensor  		&  $\tenpi_\S^\Cat = 0.09819225$ & $\uunits{\Uone}$    \\[0.1cm]
	cathode interfacial area factor 	&  $\theta_\iAE^\Cat = 1.96328590$	& $\uunits{\Uone}$    \\[0.1cm]
	\hline		
	separator electrolyte phase fraction  &  $\psi_\E^\Sep = 0.72713951$			& $\uunits{\Uone}$   \\[0.1cm]
	separator electrolyte porosity tensor &  $\tenpi_\E^\Sep = 0.86842790$			& $\uunits{\Uone}$   \\[0.1cm]
	\hline
	anode unit cell width 				&  $\ell^\An  = 10\in 100$ & $\uunits{\Unm}$  \\
	anode particle radius 			&  $\tilde r_\A^\Cat = 0.4$ & $\uunits{\Uone}$ \\
	anode electrolyte phase fraction  	&  $\psi_\E^\An = 0.72713951$  & $\uunits{\Uone}$ 	 \\
	anode electrolyte porosity tensor 	&  $\tenpi_\E^\An  = 0.86842790 $ & $\uunits{\Uone}$    \\[0.1cm]
	anode solid phase fraction   		&  $\psi_\S^\An = 0.27286022$ & $\uunits{\Uone}$   \\[0.1cm]
	anode solid porosity tensor  		&  $\tenpi_\S^\An = 0.09819225$ & $\uunits{\Uone}$      \\[0.1cm]
	anode interfacial area factor 		&  $\theta_\iAE^\An = 1.96328590$  & $\uunits{\Uone}$  	 \\[0.1cm]
\hline
\end{longtable}


\bibliographystyle{siamplain}

\begin{thebibliography}{10}

\bibitem{ngsolve}
NGSolve - Finite elemente software \url{https://ngsolve.org/}.

\bibitem{pymor}
pyMOR - Model order reduction with Python. \url{http://www.pymor.org}.

\bibitem{doi:10.1137/0523084}
G.~Allaire.
\newblock Homogenization and two-scale convergence.
\newblock {\em SIAM Journal on Mathematical Analysis}, 23(6):1482--1518, 1992.

\bibitem{eim}
M.~Barrault, Y.~Maday, N.~Nguyen, and A.~Patera.
\newblock An `empirical interpolation' method: application to efficient
  reduced-basis discretization of partial differential equations.
\newblock {\em Comptes Rendus Mathematique}, 339:667--672, 2004.

\bibitem{BARRE2013680}
A.~Barr{\'e}, B.~Deguilhem, S.~Grolleau, M.~G{\'e}rard, F.~Suard, and D.~Riu.
\newblock A review on lithium-ion battery ageing mechanisms and estimations for
  automotive applications.
\newblock {\em Journal of Power Sources}, 241:680--689, 2013.

\bibitem{doi:10.1021/ar300145c}
M.~Z. Bazant.
\newblock Theory of chemical kinetics and charge transfer based on
  nonequilibrium thermodynamics.
\newblock {\em Accounts of Chemical Research}, 46(5):1144--1160, 2013.
\newblock PMID: 23520980.

\bibitem{MR3672144}
P.~Benner, A.~Cohen, M.~Ohlberger, and K.~Willcox, editors.
\newblock {\em Model reduction and approximation}, volume~15 of {\em
  Computational Science \& Engineering}.
\newblock Society for Industrial and Applied Mathematics (SIAM), Philadelphia,
  PA, 2017.
\newblock Theory and algorithms.

\bibitem{MR3701994}
P.~Benner, M.~Ohlberger, A.~Patera, G.~Rozza, and K.~Urban, editors.
\newblock {\em Model reduction of parametrized systems}, volume~17 of {\em
  MS\&A. Modeling, Simulation and Applications}.
\newblock Springer, Cham, 2017.
\newblock Selected papers from the 3rd MoRePaS Conference held at the
  International School for Advanced Studies (SISSA), Trieste, October 13--16,
  2015.

\bibitem{HAPOD}
S.~Rave C.~Himpe, T.~Leibner.
\newblock Hierarchical approximate proper orthogonal decomposition.
\newblock {\em SIAM J. Sci. Comput.}, 40(5):3267--3292, 2018.

\bibitem{Cahn:1959ab}
J.W. Cahn.
\newblock Free energy of a nonuniform system. ii. thermodynamic basis.
\newblock {\em J. Chem. Phys.}, 30(5):1121--1124, 1959.

\bibitem{Cahn:1958aa}
J.W. Cahn and J.E. Hilliard.
\newblock Free energy of a nonuniform system. i. interfacial free energy.
\newblock {\em J. Chem. Phys.}, 28(2):258--267, 1958.

\bibitem{White}
L.~Cai and R.~White.
\newblock Reduction of model order based on proper orthogonal decomposition for
  lithium-ion battery simulations.
\newblock {\em Journal of The Electrochemical Society}, 156(3):154--161, 2009.

\bibitem{deim}
S.~Chaturantabut and D.~Sorensen.
\newblock Nonlinear model reduction via discrete empirical interpolation.
\newblock {\em SIAM J. Sci. Comput.}, 32(5):2737--2764, 2010.

\bibitem{Chung_2013}
D.-W. Chung, M.~Ebner, D.R. Ely, V.~Wood, and R.~E. Garc{\'{\i}}a.
\newblock Validity of the bruggeman relation for porous electrodes.
\newblock {\em Modelling and Simulation in Materials Science and Engineering},
  21(7):074009, 2013.

\bibitem{dGM84}
S.R. de~Groot and P.~Mazur.
\newblock {\em Non-Equilibrium Thermodynamics}.
\newblock {Dover Publications}, 1984.

\bibitem{DEFAY1977498}
R~Defay, I~Prigogine, and A~Sanfeld.
\newblock Surface thermodynamics.
\newblock {\em Journal of Colloid and Interface Science}, 58(3):498 -- 510,
  1977.

\bibitem{Newman}
M.~Doyle, T.~Fuller, and J.~Newman.
\newblock Modeling of galvanostatic charge and discharge of the
  lithium/polymer/insertion cell.
\newblock {\em Journal of the Electrochemical Society}, 140(6):1526--1533,
  1993.

\bibitem{DREYER:2011aa}
W.~Dreyer, M.~Gaberescek, C.~Guhlke, R.~Huth, and J.~Jamnik.
\newblock Phase transition in a rechargeable lithium battery.
\newblock {\em European J. Appl. Math.}, 22(03):267--290, 2011.

\bibitem{dgm2018}
W.~Dreyer, C.~Guhke, M.~Landstorfer, and R.~M{\"u}ller.
\newblock New insights on the interfacial tension of electrochemical interfaces
  and the lippmann equation.
\newblock {\em European Journal of Applied Mathematics}, 29(4):708--753, 2018.

\bibitem{Dreyer:2011ab}
W.~Dreyer, C.~Guhlke, and R.~Huth.
\newblock The behavior of a many-particle electrode in a lithium-ion battery.
\newblock {\em Physica D}, 240(12):1008 -- 1019, 2011.

\bibitem{Dreyer:2014fk}
W.~Dreyer, C.~Guhlke, and M.~Landstorfer.
\newblock A mixture theory of electrolytes containing solvation effects.
\newblock {\em Electrochemistry Communications}, 43(0):75 -- 78, 2014.

\bibitem{C3CP44390F}
W.~Dreyer, C.~Guhlke, and R.. Muller.
\newblock Overcoming the shortcomings of the nernst-planck model.
\newblock {\em Phys. Chem. Chem. Phys.}, 15:7075--7086, 2013.

\bibitem{C5CP03836G}
W.~Dreyer, C.~Guhlke, and R.~Muller.
\newblock Modeling of electrochemical double layers in thermodynamic
  non-equilibrium.
\newblock {\em Phys. Chem. Chem. Phys.}, pages~--, 2015.

\bibitem{C6CP04142F}
W.~Dreyer, C.~Guhlke, and R.~Muller.
\newblock A new perspective on the electron transfer: recovering the
  butler-volmer equation in non-equilibrium thermodynamics.
\newblock {\em Phys. Chem. Chem. Phys.}, 18:24966--24983, 2016.

\bibitem{Dreyer:2018aa}
W.~Dreyer, C.~Guhlke, and R.~M{\"u}ller.
\newblock Bulk-surface electro-thermodynamics and applications to
  electrochemistry.
\newblock {\em WIAS Preprint 2511}, 2018.

\bibitem{eim_operator}
M.~Drohmann, B.~Haasdonk, and M.~Ohlberger.
\newblock Reduced basis approximation for nonlinear parametrized evolution
  equations based on empirical operator interpolation.
\newblock {\em SIAM J. Sci. Comput.}, 34(2):937--969, 2012.

\bibitem{Ebner_2014}
M.~Ebner and V.~Wood.
\newblock Tool for tortuosity estimation in lithium ion battery porous
  electrodes.
\newblock {\em Journal of The Electrochemical Society}, 162(2):A3064--A3070,
  2014.

\bibitem{D1CP00359C}
J.S. Edge, S.~O'Kane, R.~Prosser, N.D. Kirkaldy, A.N. Patel, A.~Hales,
  A.~Ghosh, W.~Ai, J.~Chen, J.~Yang, S.~Li, M.-C. Pang, L.~Bravo~Diaz,
  A.~Tomaszewska, M.W. Marzook, K.N. Radhakrishnan, H.~Wang, Y.~Patel, B.~Wu,
  and G.J. Offer.
\newblock Lithium ion battery degradation: what you need to know.
\newblock {\em Phys. Chem. Chem. Phys.}, 23:8200--8221, 2021.

\bibitem{multibat}
J.~Feinauer, S.~Hein, S.~Rave, S.~Schmidt, D.~Westhoff, J.~Zausch, O.~Iliev,
  A.~Latz, M.~Ohlberger, and V.~Schmidt.
\newblock {MULTIBAT}: Unified workflow for fast electrochemical {3D}
  simulations of lithium-ion cells combining virtual stochastic
  microstructures, electrochemical degradation models and model order
  reduction.
\newblock {\em Journal of Computational Science}, 31:172--184, 2019.

\bibitem{Fuoss1978}
R.~M. Fuoss.
\newblock Review of the theory of electrolytic conductance.
\newblock {\em Journal of Solution Chemistry}, 7(10):771--782, Oct 1978.

\bibitem{doi:10.1021/j150551a038}
R.M. Fuoss and L.~Onsager.
\newblock Conductance of unassociated electrolytes.
\newblock {\em The Journal of Physical Chemistry}, 61(5):668--682, 1957.

\bibitem{Garcia:2004aa}
R.E. Garcia, C.M. Bishop, and W.C. Carter.
\newblock Thermodynamically consistent variational principles with applications
  to electrically and magnetically active systems.
\newblock {\em Acta Mater.}, 52(1):11--21, January 2004.

\bibitem{Pod_Greedy_RE}
B.~Haasdonk.
\newblock Convergence rates of the pod-greedy method.
\newblock {\em M2AN Mathematical Modelling and Numerical Analysis},
  47:859--873, 2013.

\bibitem{RB_Haasdonk}
B.~Haasdonk.
\newblock Reduced basis methods for parametrized pdes---a tutorial introduction
  for stationary and instationary problems.
\newblock In P.~Benner, M.~Ohlberger, A.~Cohen, and K.~Willcox, editors, {\em
  Model Reduction and Approximation}, chapter~2, pages 65--136. SIAM
  Computational Science and Engineering, 2017.

\bibitem{Pod_Greedy}
B.~Haasdonk and M.~Ohlberger.
\newblock Reduced basis method for finite volume approximations of parametrized
  linear evolution equations.
\newblock {\em M2AN Mathematical Modelling and Numerical Analysis},
  42(2):277--302, 2008.

\bibitem{HAN2019100005}
X.~Han, L.~Lu, Y.~Zheng, X.~Feng, Z.~Li, J.~Li, and M.~Ouyang.
\newblock A review on the key issues of the lithium ion battery degradation
  among the whole life cycle.
\newblock {\em eTransportation}, 1:100005, 2019.

\bibitem{RB_Rozza}
J.~S. Hesthaven, G.~Rozza, and B.~Stamm.
\newblock {\em Certified reduced basis methods for parametrized partial
  differential equations}.
\newblock SpringerBriefs in Mathematics. Springer International Publishing,
  2016.

\bibitem{Zhang_2012}
O.~Iliev, A.~Latz, J.~Zausch, and S.~Zhang.
\newblock On some model reduction approaches for simulations of processes in
  li-ion battery.
\newblock pages 161--171. Proceedings of Algoritmy 2012, conference on
  scientific computing, Vysok\'{e} Tatry, Podbansk\'{e}, Slovakia, Slovak
  University of Technology in Bratislava, 2012.

\bibitem{Newman:1973aa}
{J. Newman, K. Thomas}.
\newblock {\em Electrochemical Systems}.
\newblock John Wiley \& Sons, 2014.

\bibitem{Kim01012016}
S.~U. Kim and V.~Srinivasan.
\newblock A method for estimating transport properties of concentrated
  electrolytes from self-diffusion data.
\newblock {\em Journal of The Electrochemical Society}, 163(14):A2977--A2980,
  2016.

\bibitem{LIU2014447}
Jing L. and C.~W. Monroe.
\newblock Solute-volume effects in electrolyte transport.
\newblock {\em Electrochimica Acta}, 135:447 -- 460, 2014.

\bibitem{Landstorfer01012017}
M.~Landstorfer.
\newblock Boundary conditions for electrochemical interfaces.
\newblock {\em Journal of The Electrochemical Society}, 164(11):E3671--E3685,
  2017.

\bibitem{LANDSTORFER201856}
M.~Landstorfer.
\newblock On the dissociation degree of ionic solutions considering solvation
  effects.
\newblock {\em Electrochemistry Communications}, 92:56 -- 59, 2018.

\bibitem{Landstorfer_2019}
M.~Landstorfer.
\newblock A discussion of the cell voltage during discharge of an intercalation
  electrode for various c-rates based on non-equilibrium thermodynamics and
  numerical simulations.
\newblock {\em Journal of The Electrochemical Society}, 167(1):013518, nov
  2019.

\bibitem{M.-Landstorfer:2011aa}
M.~Landstorfer, S.~Funken, and T.~Jacob.
\newblock An advanced model framework for solid electrolyte intercalation
  batteries.
\newblock {\em Phys. Chem. Chem. Phys.}, 13:12817--12825, 2011.

\bibitem{Landstorfer2016187}
M.~Landstorfer, C.~Guhlke, and W.~Dreyer.
\newblock Theory and structure of the metal-electrolyte interface incorporating
  adsorption and solvation effects.
\newblock {\em Electrochimica Acta}, 201:187 -- 219, 2016.

\bibitem{Landstorfer:2013ly}
M.~Landstorfer and T.~Jacob.
\newblock Mathematical modeling of intercalation batteries at the cell level
  and beyond.
\newblock {\em Chem. Soc. Rev.}, 42:3234--3252, 2013.

\bibitem{LANDSTORFER2021110071}
M.~Landstorfer, B.~Prifling, and V.~Schmidt.
\newblock Mesh generation for periodic 3d microstructure models and computation
  of effective properties.
\newblock {\em Journal of Computational Physics}, 431:110071, 2021.

\bibitem{Lass_Volkwein}
O.~Lass and S.~Volkwein.
\newblock Parameter identification for nonlinear elliptic-parabolic systems
  with application in lithium-ion battery modeling.
\newblock {\em Computational Optimization and Applications}, 62(1):217--239,
  2015.

\bibitem{Latz:2011aa}
A.~Latz and J.~Zausch.
\newblock Thermodynamic consistent transport theory of li-ion batteries.
\newblock {\em J. Power Sources}, 196(6):3296--3302, 2011.
\newblock cited By (since 1996) 0.

\bibitem{eim_gen_operator}
Y.~Maday and O.~Mula.
\newblock A generalized empirical interpolation method: application of reduced
  basis techniques to data assimilation.
\newblock In F.~Brezzi, P.~Colli Franzone, U.~Gianazza, and G.~Gilardi,
  editors, {\em Analysis and Numerics of Partial Differential Equations},
  volume~4 of {\em Springer INdAM Series}, pages 221--235. Springer, 2013.

\bibitem{pymor2}
R.~Milk, S.~Rave, and F.~Schindler.
\newblock {pyMOR} -- generic alorithms and interfaces for model order
  reduction.
\newblock {\em SIAM J. Sci. Comput.}, 38(5):194--216, 2016.

\bibitem{Muller:1985fk}
I.~M{\"u}ller.
\newblock {\em Thermodynamics}.
\newblock Pitman, 1985.

\bibitem{Newman:1965aa}
J.~Newman.
\newblock The polarized diffuse double layer.
\newblock {\em Transactions of the Faraday Society}, 1965.

\bibitem{doi:10.1002/bbpc.19650690712}
J.~Newman, D.~Bennion, and C.~W. Tobias.
\newblock Mass transfer in concentrated binary electrolytes.
\newblock {\em Berichte der Bunsengesellschaft f{\"u}r physikalische Chemie},
  69(7):608--612, 1965.

\bibitem{Newman:1973ab}
J.~Newman and T.W. Chapman.
\newblock Restricted diffusion in binary solutions.
\newblock {\em AIChE Journal}, 19(2):343--348, 1973.

\bibitem{Stephan}
M.~Ohlberger and S.~Rave.
\newblock Localized reduced basis approximation of a nonlinear finite volume
  battery model with resolved electrode geometry.
\newblock In P.~Benner, M.~Ohlberger, A.~Patera, and K.~Urban G.~Rozza,
  editors, {\em Model Reduction of Parametrized Systems}, number~17 in MS\&A,
  pages 201--212. Springer International Publishing, 2017.

\bibitem{Stephan_Felix}
M.~Ohlberger, S.~Rave, and F.~Schindler.
\newblock Model reduction for multiscale lithium-ion battery simulation.
\newblock In B.~Karas{\"o}zen, M.~Manguoglu, M.~Tezer-Sezgin, S.~G\"{o}ktepe,
  and \"{O}m\"{u}r Ugur, editors, {\em Numerical Mathematics and Advanced
  Applications ENUMATH 2015}, volume 112 of {\em Lecture Notes in Computational
  Science and Engineering}, pages 317--331. Springer, 2016.

\bibitem{Stephan_Zhang}
M.~Ohlberger, S.~Rave, S.~Schmidt, and S.~Zhang.
\newblock A model reduction framework for efficient simulation of li-ion
  batteries.
\newblock In J.~Fuhrmann, M.~Ohlberger, and C.~Rohde, editors, {\em Finite
  Volumes for Complex Applications VII-Elliptic, Parabolic and Hyperbolic
  Problems}, volume~78 of {\em Springer Proceedings in Mathematics and
  Statistics}, pages 695--702. Springer International Publishing, 2014.

\bibitem{PELLETIER2017158}
S.~Pelletier, O.~Jabali, G.~Laporte, and M.~Veneroni.
\newblock Battery degradation and behaviour for electric vehicles: Review and
  numerical analyses of several models.
\newblock {\em Transportation Research Part B: Methodological}, 103:158--187,
  2017.
\newblock Green Urban Transportation.

\bibitem{RB_Manzoni}
A.~Quarteroni, A.~Manzoni, and F.~Negri.
\newblock {\em Reduced basis methods for partial differential equations}.
\newblock La Matematica per il 3+2. Springer International Publishing, 2016.

\bibitem{sanfeld1968introduction}
A.~Sanfeld.
\newblock {\em Introduction to the thermodynamics of charged and polarized
  layers}.
\newblock Monographs in statistical physics and thermodynamics. Wiley, 1968.

\bibitem{SANFELD1989L521}
A.~Sanfeld, A.~Passerone, E.~Ricci, and J.C. Joud.
\newblock Capillary properties and chemical reactivity: A thermodynamic study.
\newblock {\em Surface Science Letters}, 219(1):L521 -- L526, 1989.

\bibitem{POD}
L.~Sirovich.
\newblock Turbulence and the dynamics of coherent structures part i: coherent
  structures.
\newblock {\em Quarterly of Applied Mathematics}, 45(3):561--571, 1987.

\bibitem{FE}
V.~Thom{\'e}e.
\newblock {\em Galerkin finite elemente methods for parabolic problems}.
\newblock Springer Series in Computational Mathematics. Springer, 2006.

\bibitem{Tjaden:2016wk}
B.~Tjaden, S.J. Cooper, D.JL Brett, D.~Kramer, and P.R. Shearing.
\newblock On the origin and application of the bruggeman correlation for
  analysing transport phenomena in electrochemical systems.
\newblock {\em Current Opinion in Chemical Engineering}, 12:44 -- 51, 2016.
\newblock Nanotechnology / Separation Engineering.

\bibitem{Wesche_Volkwein}
A.~Wesche and S.~Volkwein.
\newblock The reduced basis method applied to transport equations of a
  lithium-ion battery.
\newblock {\em COMPEL: The International Journal for Computation and
  Mathematics in Electrical and Electronic Engineering}, 32:1760--1772, 2013.

\bibitem{Xia}
L.~Xia, E.~Najafi, Z.~Li, H.J. Bergveld, and M.C.F. Donkers.
\newblock A computationally efficient implementation of a full and
  reduced-order electrochemistry-based model for li-ion batteries.
\newblock {\em Applied Energy}, 208:1285--1296, 2017.

\end{thebibliography}

\end{document}